\PassOptionsToPackage{lowtilde}{url} 
\documentclass[a4paper]{siamart}

\usepackage[utf8]{inputenx}      
\usepackage[T1]{fontenc}         
\usepackage{graphicx}
\usepackage{amsfonts}
\usepackage{mathtools}
\usepackage{grffile}             
\usepackage{siunitx}
\usepackage{algpseudocode}
\usepackage{listings}
\usepackage{booktabs}
\usepackage{multirow}
\usepackage{enumitem}
\Crefname{enumi}{}{}

\usepackage{xcolor}
\definecolor{darkgreen}{rgb}{0.0,0.6,0.0}
\usepackage[draft]{changes}
\definechangesauthor[name={Roland Herzog},color=red]{RH}
\definechangesauthor[name={Tommy Etling},color=darkgreen]{TE}
\definechangesauthor[name={Martin Siebenborn},color=violet]{MS}
\setauthormarkup{}
\newcommand{\stkout}[1]{\ifmmode\text{\sout{\ensuremath{#1}}}\else\sout{#1}\fi}
\setdeletedmarkup{\stkout{#1}}

\usepackage{Local_Definitions}
\newboolean{ispreprint}
\setboolean{ispreprint}{false}

\newsiamthm{example}{Example}
\newsiamthm{remark}{Remark}
\newsiamthm{assumption}{Assumption}

\renewcommand{\theequation}{\thesection.\arabic{equation}} 
\renewcommand{\thefigure}{\thesection.\arabic{figure}}
\renewcommand{\thetable}{\thesection.\arabic{table}}
\numberwithin{theorem}{section}

\newcommand{\TheTitle}{Optimum Experimental Design for Interface Identification Problems} 
\newcommand{\TheAuthors}{T. Etling and R. Herzog and M. Siebenborn}

\headers{OED FOR INTERFACE IDENTIFICATION PROBLEMS}{\TheAuthors}

\title{{\TheTitle}\thanks{This version dated \today.}}

\author{
	Tommy Etling\thanks{Technische Universität Chemnitz, Faculty of Mathematics, Professorship Numerical Mathematics (Partial Differential Equations), D--09107 Chemnitz, Germany (\email{tommy.etling@mathematik.tu-chemnitz.de}, \url{https://www.tu-chemnitz.de/mathematik/part_dgl/people/etling}).}
	\and
	Roland Herzog\thanks{Technische Universität Chemnitz, Faculty of Mathematics, Professorship Numerical Mathematics (Partial Differential Equations), D--09107 Chemnitz, Germany (\email{roland.herzog@mathematik.tu-chemnitz.de}, \url{https://www.tu-chemnitz.de/herzog}).}
	\and
	Martin Siebenborn\thanks{University of Hamburg, Department of Mathematics, Bundesstraße~55, D--20146 Hamburg (\email{martin.siebenborn@uni-hamburg.de}, \url{https://www.math.uni-hamburg.de/home/siebenborn}).}
}

\ifpdf
\hypersetup{
	pdftitle={\TheTitle},
	pdfauthor={\TheAuthors}
}
\fi

\begin{document}

\maketitle


\begin{abstract}
	The identification of the interface of an inclusion in a diffusion process is considered.
	This task is viewed as a parameter identification problem in which the parameter space bears the structure of a shape manifold.
	A corresponding optimum experimental design (OED) problem is formulated in which the activation pattern of an array of sensors in space and time serves as experimental condition.
	The goal is to improve the estimation precision within a certain subspace of the infinite dimensional tangent space of shape variations to the manifold, and to find those shape variations of best and worst identifiability.
	Numerical results for the OED problem obtained by a simplicial decomposition algorithm are presented.
\end{abstract}

\begin{keywords}
	optimum experimental design, 
	diffusion process,
	sensor activation,
	inclusion detection,
	interface identification
\end{keywords}

\begin{AMS}
	62K05, 
	35R30, 
	35K20, 
	49Q10, 
	90C25  
\end{AMS}

\section{Introduction}
\label{sec:Introduction}

Optimum experimental design (OED) aims at improving the setup of experiments in order to increase the precision of parameter estimation in the face of measurement errors.
Often, the model in which the parameters are to be fitted, is given by a time-dependent ordinary or partial differential equation.
Experimental conditions which are available for optimization may include initial conditions, boundary conditions, right hand side forces, as well as sensor placement and activation patterns.
We refer the reader to \cite{Ucinski2005:1,PronzatoPazman2013} for a general introduction to OED problems.

In this paper we consider a diffusion model in which the parameter to be determined is the location of an interface surrounding an inclusion with an alternate value of the diffusion coefficient in an otherwise homogeneous body.
The state of the system represents the concentration of a diffusive substance.
Measurements of the state are taken through an array of sensors.
Since these measurements are distributed, we take into account spatially correlated measurement errors.
We seek to optimize the sensor activation pattern such that the response in the observables under variations of the parameter, and thus the precision of the estimation, are maximized.
The latter is measured in terms of the Fisher information.

The novelty of the problem under consideration lies in the fact that the space of potential interface locations, or equivalently, the space of possible shapes of the inclusion, does not bear the structure of a vector space.
By contrast, we follow \cite{MichorMumford2006,Schulz2014} and treat two-dimensional shapes as the infinite dimensional manifold of smooth embeddings of the circle $\SS^1$ modulo diffeomorphisms. 
Variations of a shape are then elements of the tangent space at that point, which can be identified with the vector space of smooth normal vector fields (aka velocity fields) along the shape boundary.
This feature sets our approach apart from common formulations of OED problems discussed in the literature, in which the underlying parameter estimation problem seeks to determine a parameter in $\R^n$, which coincides with its tangent space.

Consequently, we need to revisit and refine the usual machinery to formulate an OED problem for our interface identification problem.
Notice that since the tangent space is infinite dimensional, it contains large subspaces which are poorly identifiable.
This fact reflects the ill-posedness of the estimation problem and it cannot be overcome by design of experiments alone.
In the actual interface estimation problem, poorly identifiable subspaces are usually suppressed by a suitable regularization term, the most common of which is boundary length parametrization.

In our OED problem we follow the paradigm of regularization by discretization and focus on a finite dimensional subspace of parameter variations which are reasonably well identifiable.
Generally speaking, this subspace contains shape variations with features on a length scale somewhat smaller than the interface length.
In practice, we use smooth bump functions distributed around the interface perimeter to span this subspace.
Then we improve on the identifiability of shape variations within this subspace by optimizing the experimental conditions, i.e., the pattern of sensor activations.

Notice that when the unknown parameter belongs to $\R^n$, the standard basis of the latter is chosen and the Fisher information matrix (FIM) is formed w.r.t.\ this basis.
The OED criterion to be optimized is then formulated in terms of the eigenvalues of the FIM.
In our problem, the fact that the velocity fields describing the subspace are generally not orthonormal, represents another, albeit minor novel feature of our work.
It implies that either an orthonormalization w.r.t.\ the chosen Riemannian metric in the tangent space has to be carried out before setting up the Fisher information matrix (FIM), or else a generalized eigenvalue problem for the FIM has to be considered instead of an ordinary one in order that the OED criterion becomes independent of the actual basis.

Finally, we mention that the derivative of the parameter-to-observable map, which is required in the formulation and solution of OED problems, maps velocity fields into distributed observations of the concentration.
These directional derivatives are also known as material derivatives and they are evaluated using shape optimization techniques.


Let us put our work into perspective.
Optimal sensor location/activation and other optimum experimental design problems for various identification problems in diffusion processes have been considered in the literature before; see for instance \cite{QureshiNgGoodwin1980,UcinskiDemetriou2004,WalshWildeyJakeman2017_preprint}, \cite[Chap.~8.4]{Ucinski2005:1}.
Various methods for the identification of inclusions in diffusion equations have been discussed, for instance, in \cite{ChapkoKressYoon1999,FruehaufGebauerScherzer2007,LesnicBin-Mohsin2012,HarbrechtTausch2013}.
However, we are not aware of any combination of the two techniques.

This paper is organized as follows.
In \cref{sec:Forward_and_sensitivity_problems} we describe the forward problem for the state variable as well as the sensitivity equation governing the material derivative of the state w.r.t.\ interface perturbations.
\Cref{sec:Accuracy_of_estimation} is devoted to the accuracy of the estimation through the formulation of a suitable Fisher information matrix (FIM) and design criterion.
In \cref{sec:OED_problem} we state the optimum experimental design problem and discuss an algorithm for its solution.
Numerical results are presented in \cref{sec:Numerical_results}.

\section{Forward and Sensitivity Problems}
\label{sec:Forward_and_sensitivity_problems}

In this section we consider the forward problem and its sensitivities w.r.t.\ to interface, i.e., parameter perturbations.

\subsection{Forward Problem}

We consider the diffusion problem on a space-time cylinder $\Omega \times (0,T)$, given by
\begin{subequations}
	\label{eq:forward_problem}
	\begin{alignat}{2}
		\dot u - \div (k \nabla u) 
		&
		= 
		0 
		& \quad & \text{in } \Omega \times (0,T),
		\label{eq:forward_problem_1}
		\\
		u
		& 
		= 
		\uD
		& & \text{on } \GammaD \times (0,T),
		\label{eq:forward_problem_2}
		\\
		\frac{\partial u}{\partial n} + \beta \, u
		& 
		= 
		0
		& & \text{on } \GammaR \times (0,T),
		\label{eq:forward_problem_3}
		\\
		u(\cdot,0)
		& 
		= 0
		& & \text{in } \Omega.
		\label{eq:forward_problem_4}
	\end{alignat}
\end{subequations}

\begin{figure}[htbp]
	\centering
	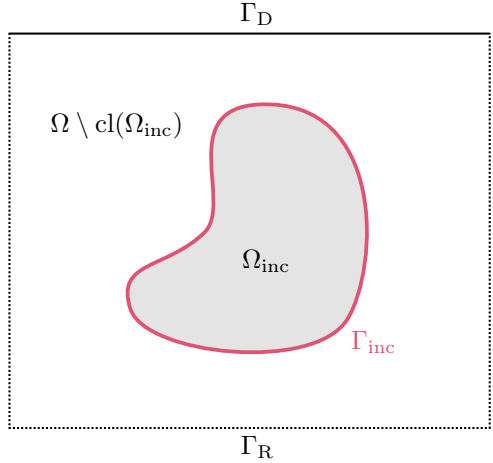
	\caption{Example for the diffusion problem described in \eqref{eq:forward_problem}.}
	\label{fig:example_domain}
\end{figure}

We work under the following assumptions on the data.
\begin{assumption} \hfill
	\label{assumption:main}
	\begin{enumerate}[label=$(\roman*)$,leftmargin=*]
		\item 
			\label{item:geometry}
			Let $\Omega \subset \R^d$, $d \in \{2,3\}$, be a bounded domain with Lipschitz boundary $\Gamma$.
			Suppose that $\Omegainc$ is an open (not necessarily connected) subset of $\Omega$ with smooth boundary $\Gammainc$ such that $\cl(\Omegainc) \subset \Omega$.
			
		\item 
			The boundary $\Gamma$ of $\Omega$ is divided into two disjoint parts $\GammaD$ and $\GammaR$, where $\GammaD$ has positive Lebesgue surface measure.

		\item 
			We assume that the Dirichlet data $\uD$ is the trace on $\Gamma_D$ of a function in $H^1(\Omega)$, also termed $\uD$.

		\item 
			The exchange coefficient $\beta$ is a non-negative function in $L^\infty(\GammaR)$.

		\item
			\label{item:diffusion_parameter}
			The diffusion parameter $k \in L^\infty(\Omega)$ is piecewise constant, i.e.,
			\begin{equation}
				\label{eq:diffusion_parameter}
				k = \chi_{\Omega\setminus\cl(\Omegainc)} \, \kbulk + \chi_{\Omegainc} \kinc,
			\end{equation}
			where $\chi_A$ denotes the characteristic function of a set $A$ and $\kbulk \neq \kinc$ are positive numbers.

	\end{enumerate}
\end{assumption}
$\Omegainc$ represents the region occupied by the inclusion.
Notice that \cref{item:geometry} implies that $\Omegainc$ has a positive distance to the boundary $\Gamma$ of $\Omega$.
We refer to \cref{fig:example_domain} for an example.

We define the spaces $\VV := H^1_D(\Omega) := \{ u \in H^1(\Omega): u = 0 \text{ on } \Gamma_D \}$ and 
\begin{equation*}
	\WW(0,T) 
	=
	\{ u \in H^1(0,T;\VV): \dot u \in L^2(0,T;\VV^*) \}
	.
\end{equation*}
The space $\WW(0,T)$ embeds continuously into $C([0,T];L^2(\Omega))$ so that the initial conditions as well as snapshots of the state at given points in time are well defined.

A weak formulation of \eqref{eq:forward_problem} is given as follows.
\begin{multline}
	\label{eq:forward_problem_weak}
	\text{Find $u \in u_D + \WW(0,T)$ such that $u(\cdot,0) = 0$ holds and}
	\\
	\dual{\dot u}{v}
	+ 
	\int_\Omega k \, (\nabla u)^\top \nabla v \, \dx
	+ 
	\int_{\GammaR} \beta \, u \, v \, \ds
	= 
	0
	\quad \text{for all } v \in \VV \text{ and a.a.\ } t \in (0,T)
	.
\end{multline}
It is well known that under the given assumptions, the forward problem \eqref{eq:forward_problem} has a unique weak solution; see for instance \cite[Chap.~IV]{Wloka1982} or \cite[Chap.~7.1]{Evans1998}.

\subsection{Shape Space as a Parameter Space}
\label{subsec:Shape_space}

As we mentioned in the introduction, the space of potential interface locations, or equivalently, the space of possible shapes of the inclusion, does not bear the structure of a vector space.
By contrast, we follow \cite{MichorMumford2006,Schulz2014} and treat (connected) shapes as an infinite dimensional manifold of smooth embeddings of the circle $\SS^1$ modulo diffeomorphisms.

Variations of such a shape are then elements of the tangent space at that point, which can be identified with the vector space of all smooth normal vector fields (aka velocity fields) along the shape boundary.
This feature sets our approach apart from previous formulations of OED problems discussed in the literature, in which the underlying parameter estimation problem seeks to determine a parameter in $\R^n$, which coincides with its tangent space.

The tangent space of normal velocity fields has to be endowed with a Riemannian metric (inner product).
As in \cite{SchulzSiebenbornWelker2015:2}, we proceed as follows.
Suppose that $D$ is a hold-all domain containing $\Omegainc$ in its interior; see \cref{fig:domain_experiment}.
We set 
\begin{equation*}
	\ZZ(\bar V)
	\coloneqq
	\{ V \in H^1(D;\R^d): V = 0 \text{ on } \partial D, \; V = \bar V \text{ on } \Gammainc \}
\end{equation*}
for $\bar V \in H^{1/2}(\Gammainc)$.
A normal velocity field $\bar V$ on $\Gammainc$ is extended into all of $D$ by solving the following linear elasticity problem:
\begin{equation}
	\label{eq:elasticity_problem}
	\text{find $V \in \ZZ(\bar V)$ such that }
	\int_D 
	\bvarepsilon(V) \dprod \sigma(\bvarepsilon(W)) \, \dx
	=
	0
	\quad \text{for all } W \in \ZZ(0)
	.
\end{equation}
Here the strain $\bvarepsilon(V) \coloneqq (1/2)(\nabla V + \nabla V^\top)$ is the symmetric part of the Jacobian of $V$, and the stress-strain relation is given by $\sigma(\bvarepsilon) \coloneqq 2 \, \mu \, \bvarepsilon + \lambda \, \trace(\bvarepsilon) \, \id$.
The Lam{\'e} parameters $(\lambda,\mu)$ for this problem will be specified in \cref{sec:Numerical_results}.

From now on, we will always consider normal velocity fields to be extended as above.
Now if $V_1$ and $V_2$ are two such fields (pertaining to $\bar V_1$ and $\bar V_2$, respectively), we utilize the following inner product:
\begin{equation}
	\label{eq:Riemannian_metric}
	\manifoldinprod{V_1}{V_2}
	\coloneqq
	\int_D \nabla V_1 \dprod \nabla V_2 \, \dx
	.
\end{equation}
Notice that this also serves as an inner product between $\bar V_1$ and $\bar V_2$ in the tangent space.
The inner product \eqref{eq:Riemannian_metric} will be required in \cref{sec:Accuracy_of_estimation} in order to define the notion of an orthonormal basis.

\subsection{Sensitivity Problem}

In this section we consider the sensitivity of the forward problem \eqref{eq:forward_problem} w.r.t.\ parameter perturbations $V \in W^{1,\infty}(\Omega;\R^d)$.
These sensitivites are also known as material derivatives.
In weak formulation, the material derivative $\delta u$ in the direction of $V$ is given as follows.
\begin{multline}
	\label{eq:sensitivity_problem_weak}
	\text{Find $\delta u \in \WW(0,T)$ such that $\delta u(\cdot,0) = 0$ and for all $v \in \VV$ and a.a.\ $t \in (0,T)$,}
	\\
	\hspace*{-30mm}
	\dual{\dot{ \delta u }}{v}
	+ 
	\int_\Omega k \, (\nabla \delta u)^\top \nabla v \, \dx
	+ 
	\int_{\GammaR} \beta \, \delta u \, v \, \ds
	\\
	= 
	- \dual{\dot u}{v \, (\div V)}
	+ \int_\Omega k \, (\nabla u)^\top \left[ DV + DV^\top - (\div V) \, \id \right] \nabla v \, \dx
	.
\end{multline}
In \eqref{eq:sensitivity_problem_weak}, $DV$ and $DV^\top$ denote the Jacobian of the vector field $V$ and its transpose.
The derivation of \eqref{eq:sensitivity_problem_weak} is rather standard, see for instance \cite[Chap.~2.27]{SokolowskiZolesio1992}.

\section{Accuracy of Estimation}
\label{sec:Accuracy_of_estimation}

We recall that it is our goal to estimate the interface position from a number of state measurements, taken at different locations and time instances.
In our setup, we consider distributed measurements in a number of subdomains which can be activated independently at a number of time points.
An individual measurement of the state $u$ in the measurement domain $\Omegaobs^k$ at time $t^\ell$ is described as
\begin{equation}
	\label{eq:individual_measurement}
	E_{k,\ell} \, u 
	\coloneqq
	\restr{u(t^\ell)}{\Omegaobs^k}
	\in 
	L^2(\Omegaobs^k)
	.
\end{equation}

\begin{figure}[htbp]
	\centering
	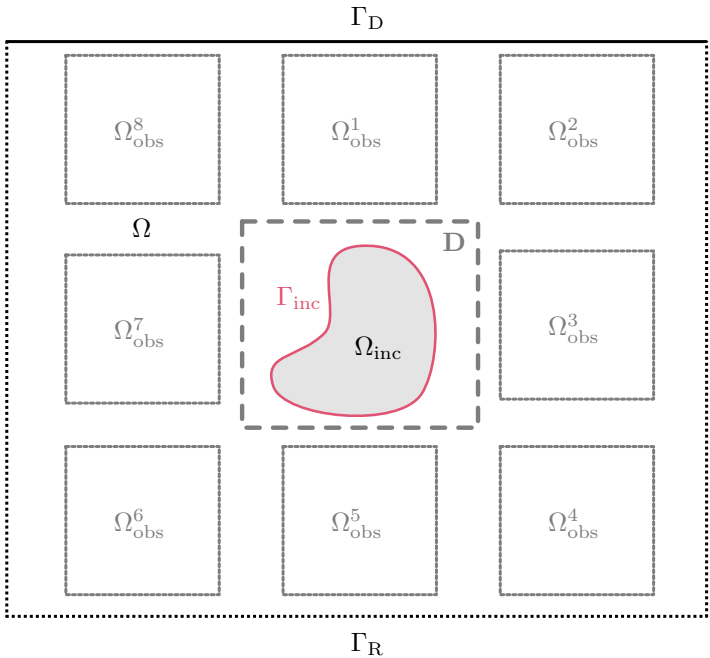
	\caption{Experimental setting with measurement domains $\Omegaobs^k$ and hold-all domain $D$.}
	\label{fig:domain_experiment}
\end{figure}

We work under the following assumptions.
\begin{assumption} \hfill
	\label{assumption:measurements}
	\begin{enumerate}[label=$(\roman*)$,leftmargin=*]
		\item 
			The spatial measurement domains $\Omegaobs^k$ are disjoint measurable subsets of $\Omega \setminus \cl(\Omegainc)$, $k = 1, \ldots, \Nobs$. 
		\item 
			The measurement times $t^\ell$ are distinct time points in $(0,T]$, $\ell = 1, \ldots, \Ntime$.
	\end{enumerate}
\end{assumption}
A typical setup is depicted in \cref{fig:domain_experiment}.

Since $E_{k,\ell}$ represents a distributed measurement, we have to anticipate spatial correlations within each measurement region.
We follow \cite{AlexanderianPetraStadlerGhattas2014} and use as a covariance operator $\CC_k = \AA_k^{-2}$ on $\Omegaobs^k$, where $\AA_k$ represents the second-order elliptic differential operator 
\begin{equation}
	\label{eq:covariance_root_differential_operator}
	\AA_k = - \alpha_0 \laplace + \alpha_1 \id
\end{equation}
defined on a dense subset of $L^2(\Omegaobs^k)$.
The positive parameters $\alpha_0$ and $\alpha_1$ encode the variation and correlation length.
On the other hand we assume that distinct measurement regions are sufficiently separated both in space and time such that no correlations occur between them.

In order to assess the quality of an experiment, one usually considers its Fisher information; see for instance \cite[Chap.~1.5--1.7]{FedorovLeonov2014:1}.
In our context, it is useful to think of the Fisher information as a bilinear form, accepting two directional derivatives of the parameter-to-observable map. 
As we explained in \cref{subsec:Shape_space}, the tangent space of the parameter shape space can be represented as the vector space of all smooth velocity fields $\bar V$ normal to the interface $\Gammainc$.
In view of the linearity of the observation operator, a directional derivative of the parameter-to-observable map has the representation
\begin{equation*}
	V \mapsto E_{k,\ell} \, \delta u,
\end{equation*}
where $V$ is the extension of $\bar V$ in virtue of \eqref{eq:elasticity_problem}, and $\delta u$ is the solution of the sensitivity equation \eqref{eq:sensitivity_problem_weak}.
The Fisher information of an elementary experiment, consisting of the single observation defined by a fixed $E_{k,\ell}$, is thus given by the bilinear form
\begin{equation}
	\label{eq:elementary_FI_operator}
	\begin{aligned}
		\overline \Upsilon_{k,\ell}(V_1,V_2)
		&
		:=
		\scalarprod{\CC_k^{-1} E_{k,\ell} \, \delta u_1}{E_{k,\ell} \, \delta u_2}_{L^2(\Omegaobs^k)}
		\\
		& 
		=
		\int_{\Omegaobs^k} (\AA_k \, \delta u_1(\cdot,t^\ell)) \, (\AA_k \, \delta u_2(\cdot,t^\ell)) \, \dx
		.
	\end{aligned}
\end{equation}
The inverse of the Fisher information operator $\overline \Upsilon_{k,\ell}$, if it exists, can be conceived as the covariance operator of the least-squares estimator based on the linearized parameter-to-observable map.
In our situation, however, the tangent space has infinite dimension.
Owing to the smoothness of $V$, and based on compact embeddings, one can show that there exist bounded sequences $\{V_j\}$ such that $\overline \Upsilon_{k,\ell}(V_j,V_j) \to 0$.
This is an expression of the ill-posedness of the interface identification problems.
In particular, high-frequency oscillations of the interface position are likely to be estimable only poorly.

\begin{figure}[htbp]
	\label{fig:untransformed_bumps}
	\begin{center}
		\begin{tabular}{ccc}
			\includegraphics[width=.3\textwidth]{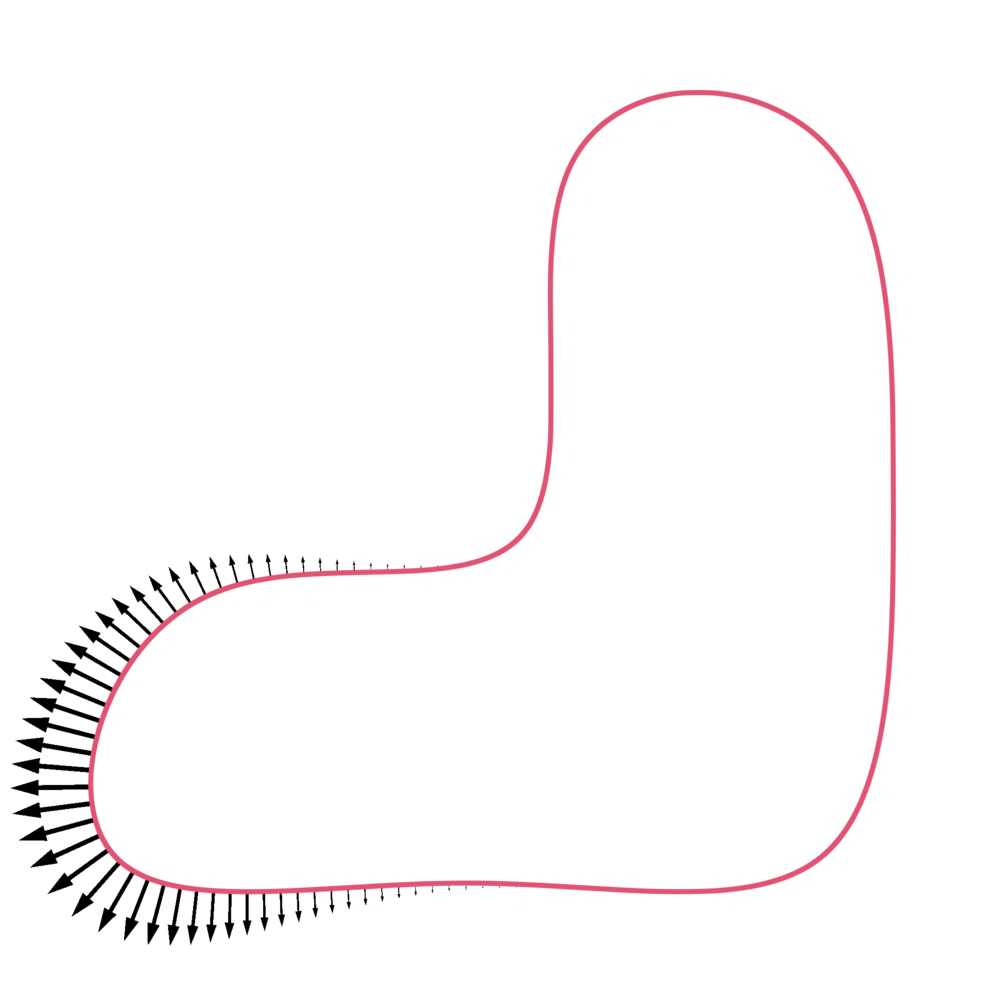}&
			\includegraphics[width=.3\textwidth]{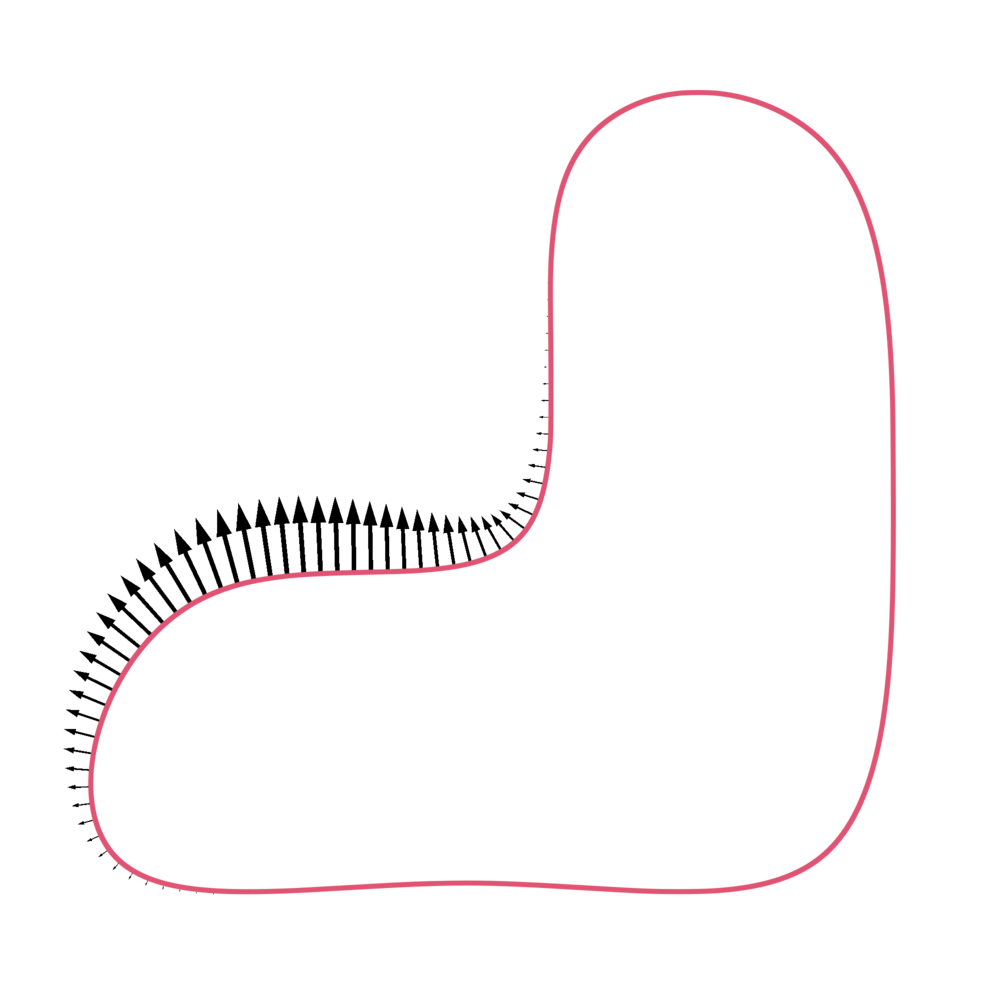}&
			\includegraphics[width=.3\textwidth]{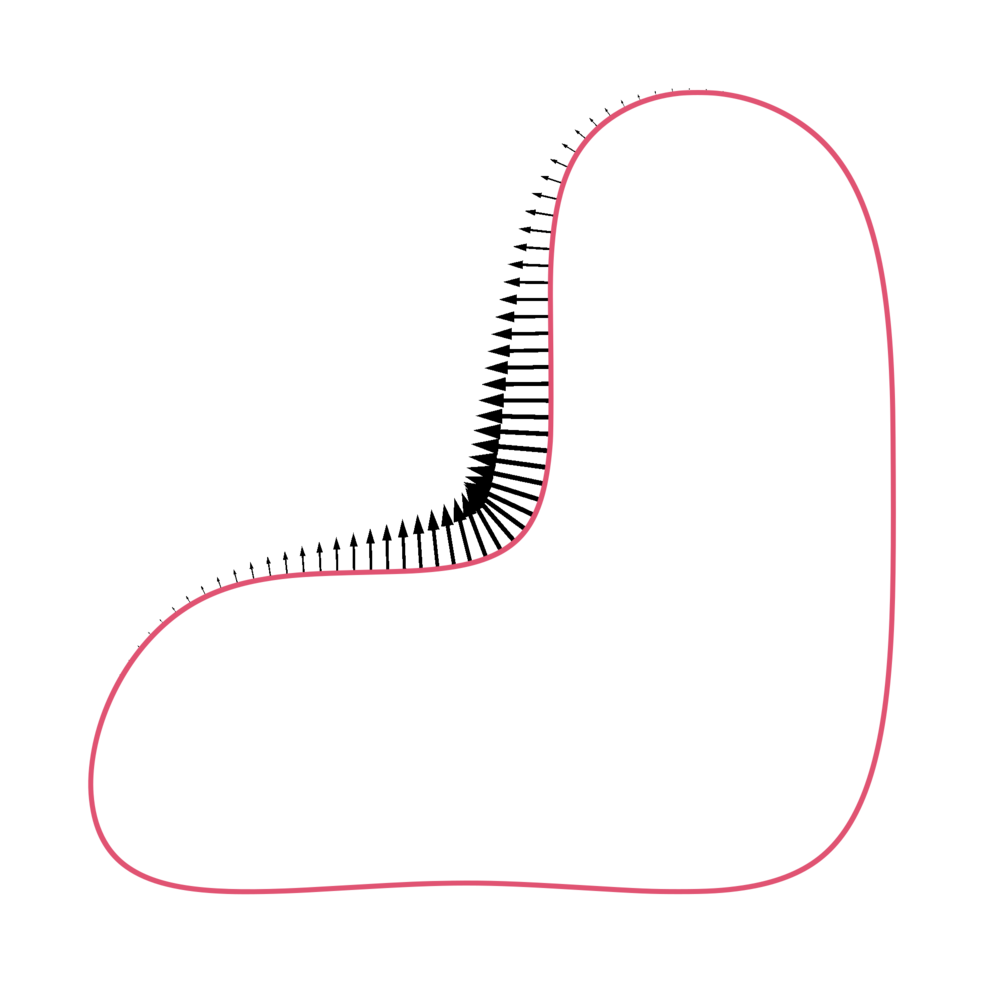}\\
			\includegraphics[width=.3\textwidth]{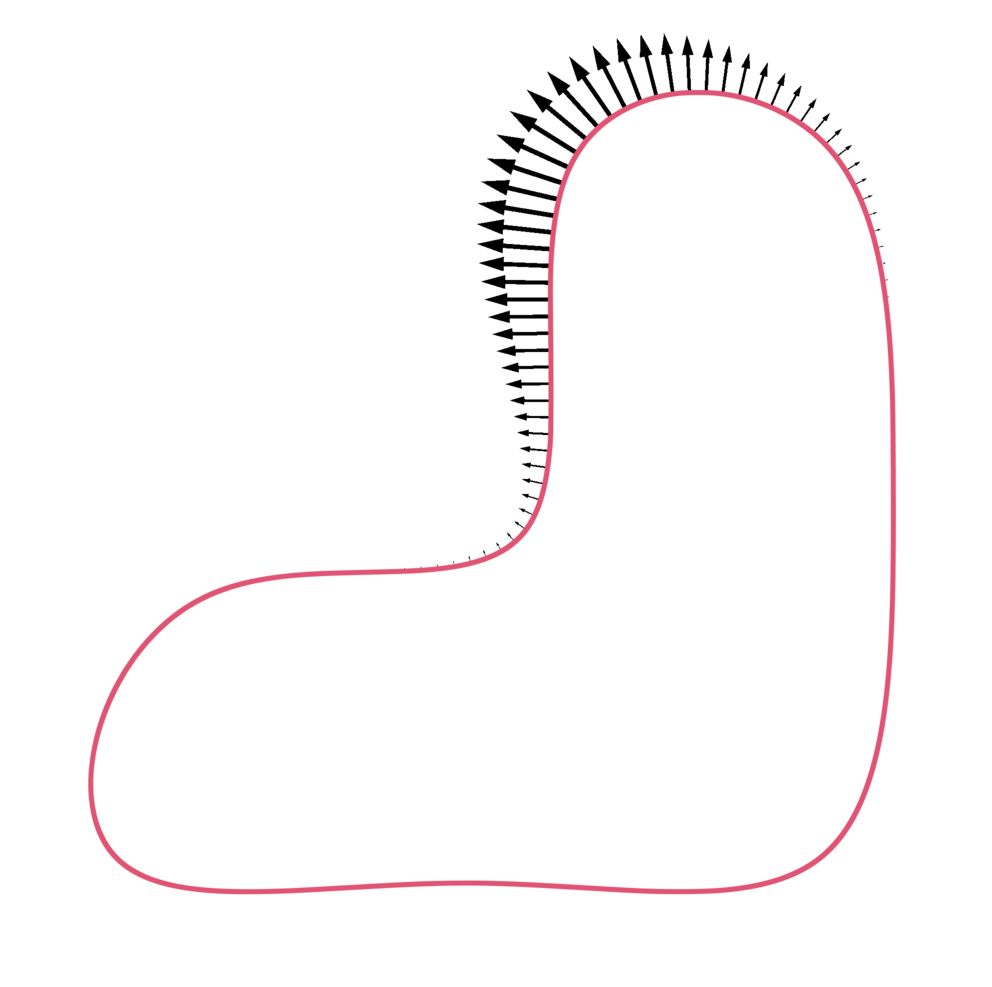}&
			\includegraphics[width=.3\textwidth]{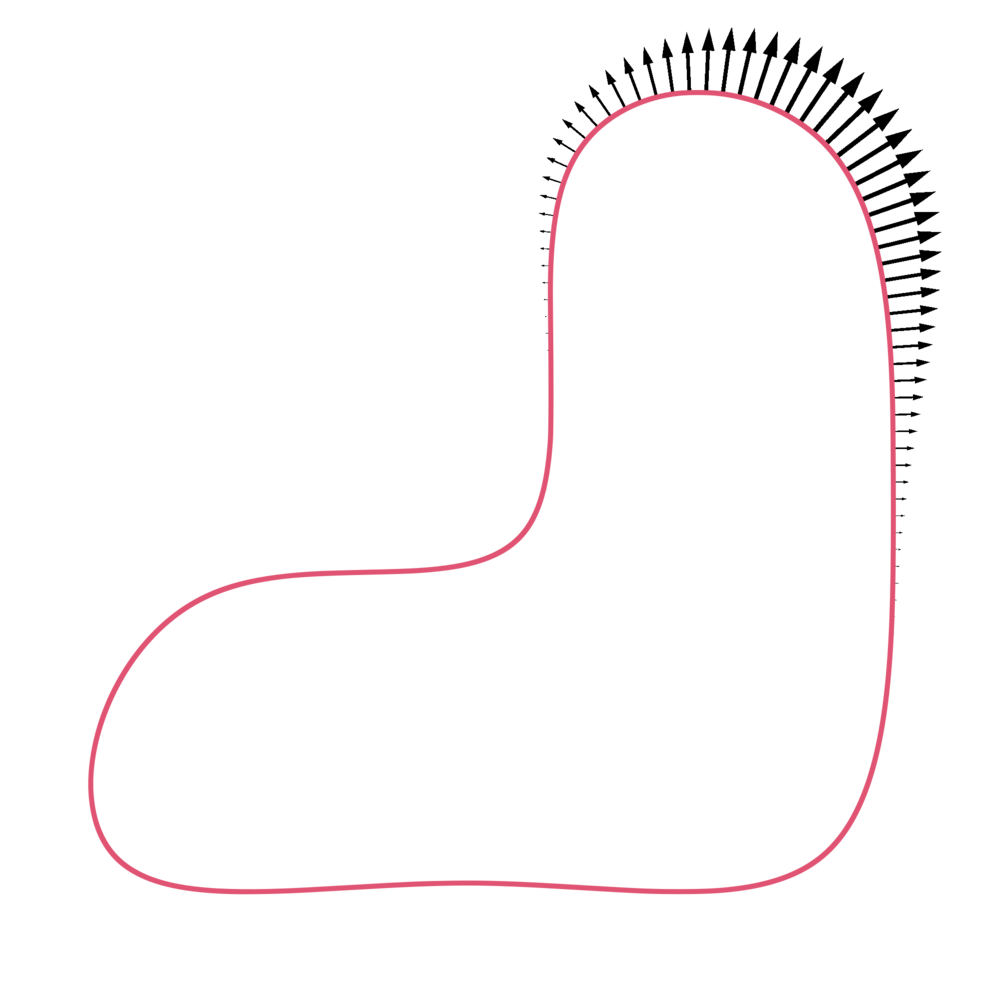}&
			\includegraphics[width=.3\textwidth]{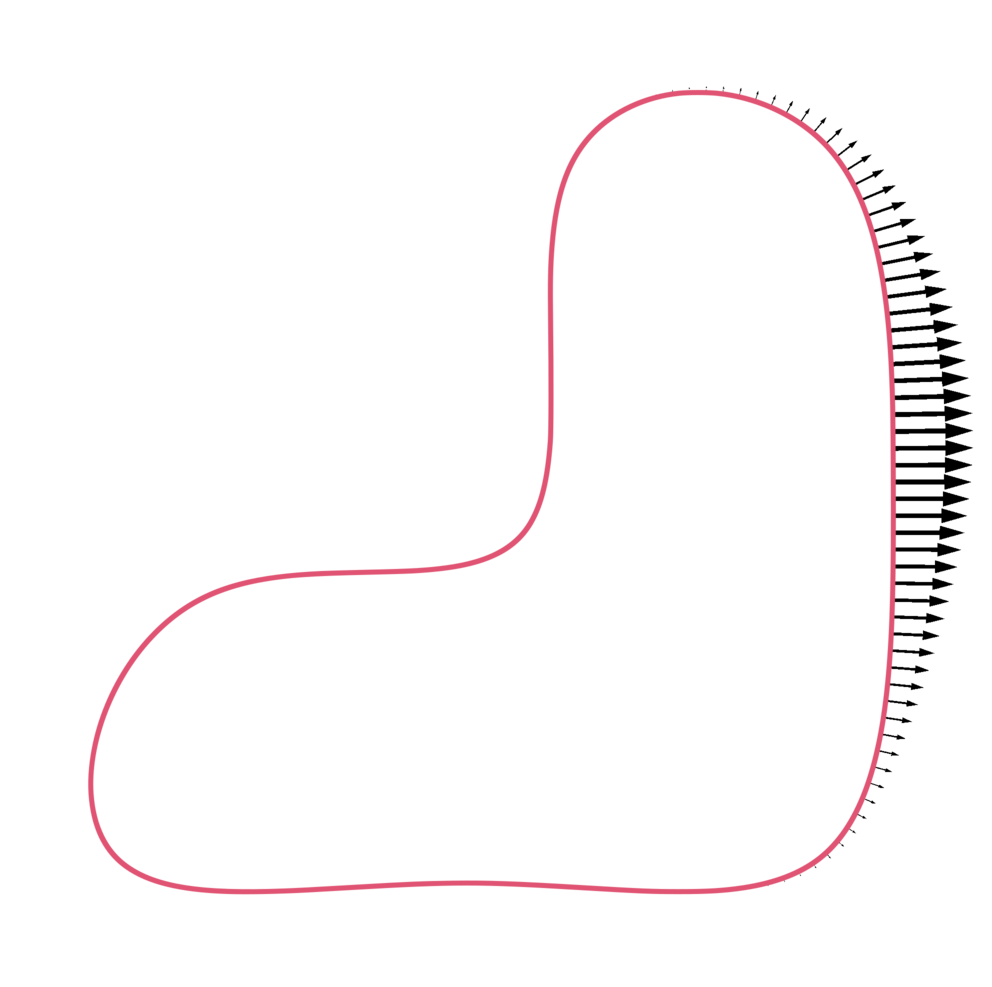}\\
			\includegraphics[width=.3\textwidth]{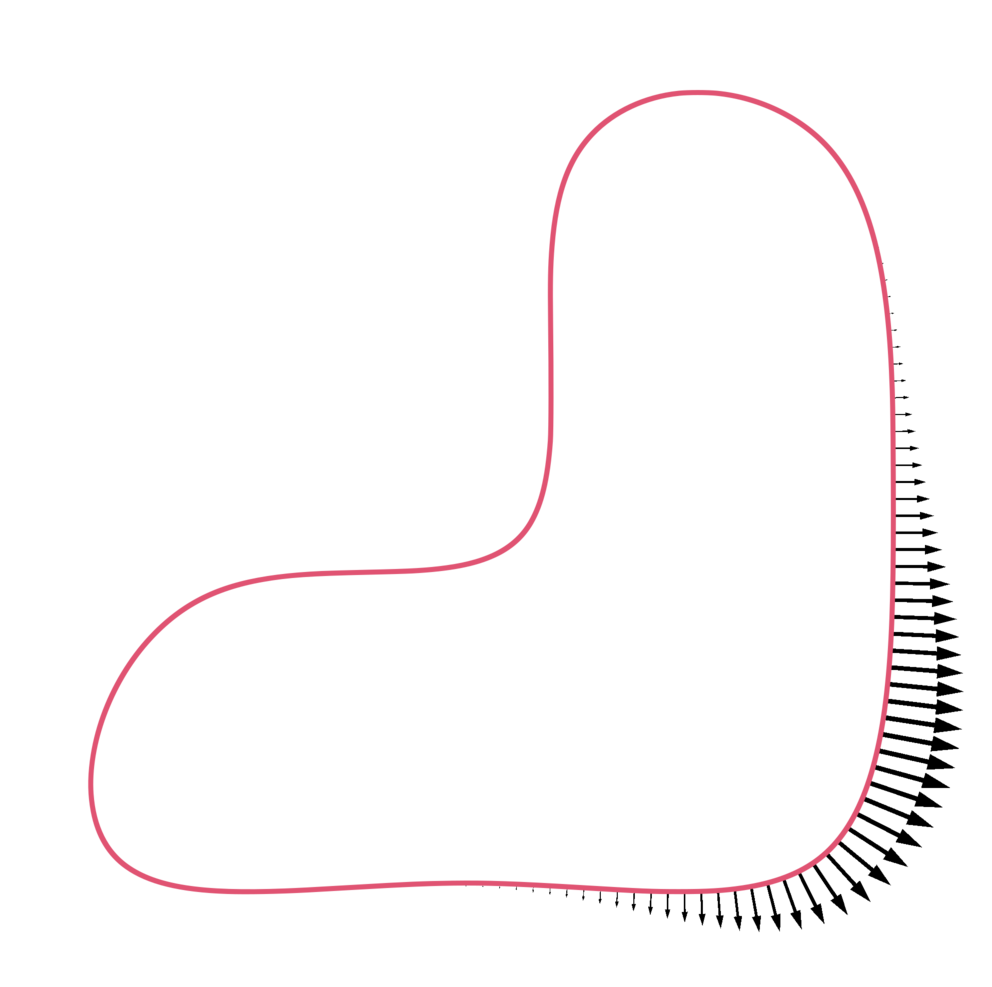}&
			\includegraphics[width=.3\textwidth]{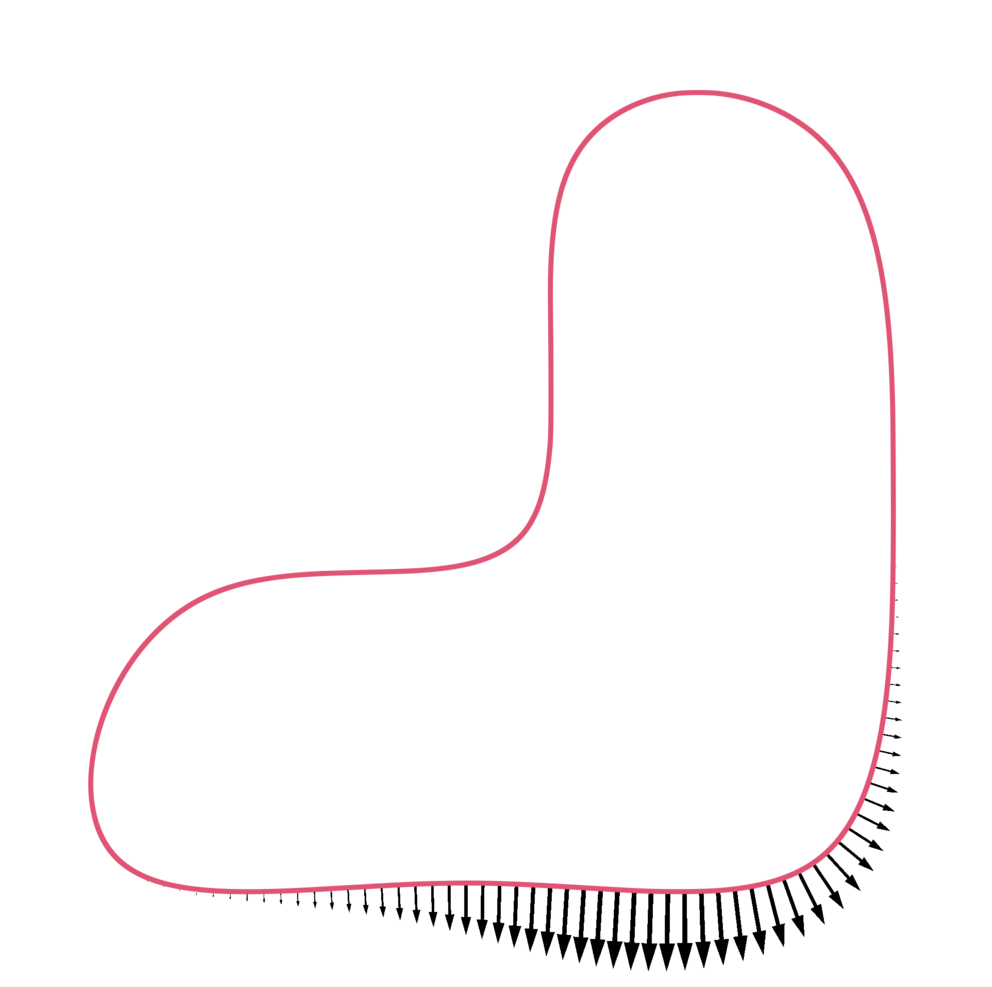}&
			\includegraphics[width=.3\textwidth]{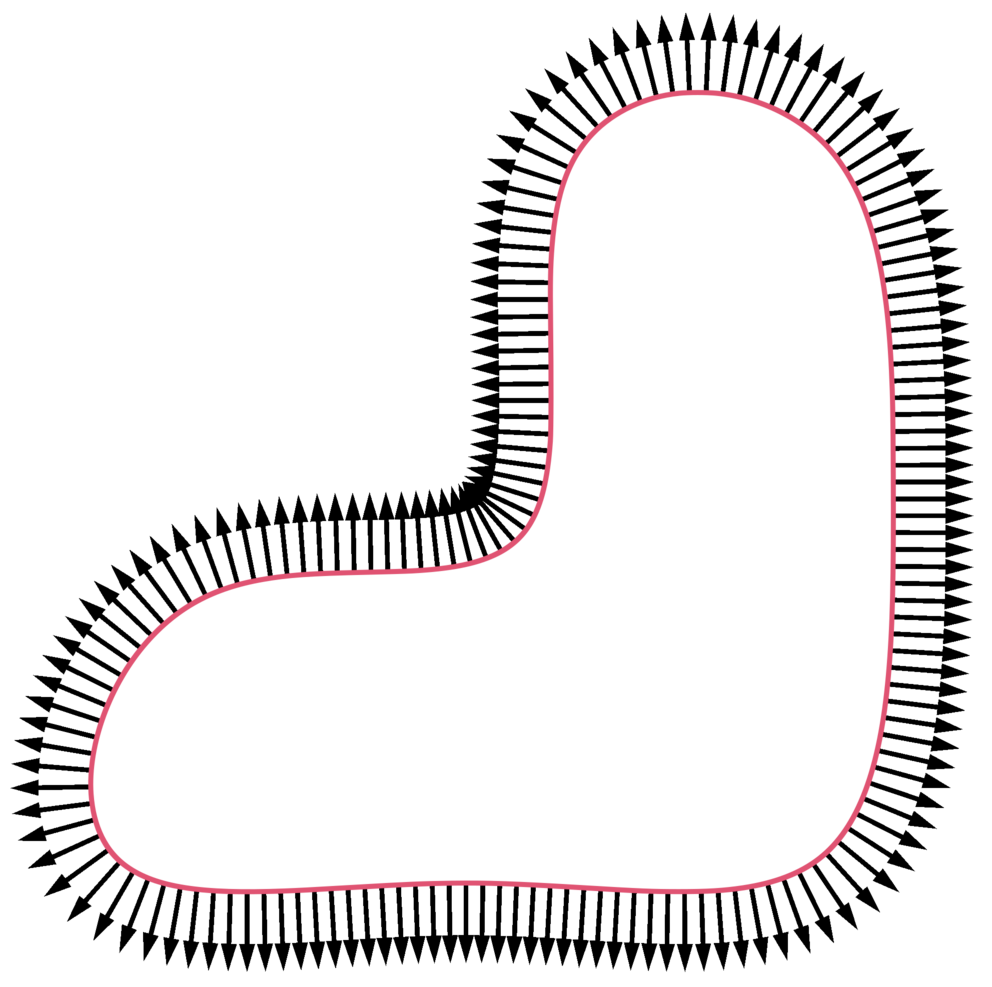}
		\end{tabular}
		\caption{Basis functions $\bar V_1, \ldots, \bar V_{\Nbasis}$ representing a subspace of shape variations of dimension $\Nbasis = 9$ on the interface $\Gammainc$.}
	\end{center}
\end{figure}

Unfortunately the structural ill-posedness of the problem cannot be cured by optimizing the measurements. 
Therefore, we follow the paradigm of regularization by discretization and restrict the discussion in the sequel to a finite dimensional subspace containing low-frequency interface variations.
Suppose that $\{V_1, \ldots, V_{\Nbasis} \}$ is a basis of such a subspace, see for example \cref{fig:untransformed_bumps}.

Choosing a finite dimensional subspace allows us to represent the Fisher information operator pertaining to the elementary experiment $E_{k,\ell}$ associated with the measurement field located at $\Omegaobs^k$ at time $t^\ell$, as an (elementary) Fisher information \emph{matrix},
\begin{equation}
	\label{eq:elementary_FIM}
	(\Upsilon_{k,\ell})_{i,j} := \overline \Upsilon_{k,\ell}(V_i,V_j).
\end{equation}
A full experiment will then be formed as a binary linear combination of elementary experiments.
Owing to the independence of measurement outcomes at individual regions and time instances, the combined Fisher information matrix (FIM) becomes
\begin{equation}
	\label{eq:combined_experiment}
	\Upsilon(w) := \sum_{k=1}^{\Nobs} \sum_{\ell=1}^{\Ntime} w_{k,\ell} \Upsilon_{k,\ell},
	\quad
	w_{k,\ell} \in \{0,1\}
	.
\end{equation}
A zero weight $w_{k,\ell} = 0$ means that no measurement will be taken within $\Omegaobs^k$ at time $t^\ell$, while $w_{k,\ell} = 1$ means that the sensor is active at the time.

Common criteria to assess the quality of an experiment consider the eigenvalues of \eqref{eq:combined_experiment}.
Since we are using in the tangent space the inner product \cref{eq:Riemannian_metric}, we in fact need to consider the generalized eigenvalue problem
\begin{equation}
	\label{eq:generalized_eigenvalue_problem}
	\Upsilon(w) \, V = \Lambda \, B \, V,
	\quad \text{where} \quad
	(B)_{ij} 
	:=
	\manifoldinprod{V_i}{V_j}
\end{equation}
is the symmetric and positive definite Gramian matrix associated with the basis $\{V_1, \ldots, V_{\Nbasis}\}$.

For concreteness, we consider in the sequel the A-criterion, i.e.,
\begin{equation}
	\label{eq:A-criterion}
	\Phi_A(\Upsilon) 
	:= 
	\sum_{i=1}^{\Nbasis} \Lambda_i^{-1}
	=
	\trace(B \Upsilon^{-1})
	.
\end{equation}
(When the symmetric and positive semi-definite matrix $\Upsilon$ has one or more zero eigenvalues, $\Phi_A(\Upsilon)$ is understood as $\infty$.)

The criterion \eqref{eq:A-criterion} allows us to compare two measurement setups.
Smaller values of $\Phi_A$ represent better experiments.
We mention that small $\Lambda_i^{-1}$ and the corresponding eigenfunctions of \eqref{eq:generalized_eigenvalue_problem} represent the subspace of shape variations with good identifiability.
Indeed, $\Lambda_i^{-1}$ can be interpreted as a squared semi-axis in a confidence ellipsoid associated with the experiment; see for instance \cite[Chap.~2.2.1]{FedorovLeonov2014:1} and \cite[Chap.~5.1.1]{PronzatoPazman2013}. 

\section{Optimum Experimental Design Problem and Algorithm}
\label{sec:OED_problem}

We recall that the optimum experimental design (OED) problem consists of selecting which measurement domains are to be active at what times in order to obtain a combined information matrix \eqref{eq:combined_experiment} with a minimal value of the objective \eqref{eq:A-criterion}.
In order to avoid facing an integer nonlinear problem of dimension $\Nobs \cdot \Ntime$, it is customary to consider a relaxation of the weights $w_{k,\ell} \in [0,1]$.
This leads to the following relaxed OED problem:
\begin{equation}
	\label{eq:relaxed_OED_problem}
	\begin{aligned}
		\text{Minimize} \quad & \Phi_A (\Upsilon(w))
		\\
		\text{where} \quad & \Upsilon(w) := \sum_{k=1}^{\Nobs} \sum_{\ell=1}^{\Ntime} w_{k,\ell} \Upsilon_{k,\ell} 
		\\
		\text{s.t.} \quad & 0 \le w_{k,\ell} \le 1 \quad \text{for all } k = 1, \ldots, \Nobs, \; \ell = 1, \ldots, \Ntime
		\\
		\text{and} \quad & \sum_{k=1}^{\Nobs} \sum_{\ell=1}^{\Ntime} w_{k,\ell} \le C_w
		.
	\end{aligned}
\end{equation}
In the interest of a concise notation, we will simply denote the double sum by $\sum_{k,\ell}$ in what follows.
Moreover, $w$ denotes the weight vector with components $w_{k,\ell}$ in some fixed enumeration.
The upper bound $C_w$ is chosen to be smaller than $\Nobs \cdot \Ntime$.

\begin{figure}[htbp]
	\label{fig:restricted_simplex}
	\centering
	\includegraphics[width=0.6\textwidth]{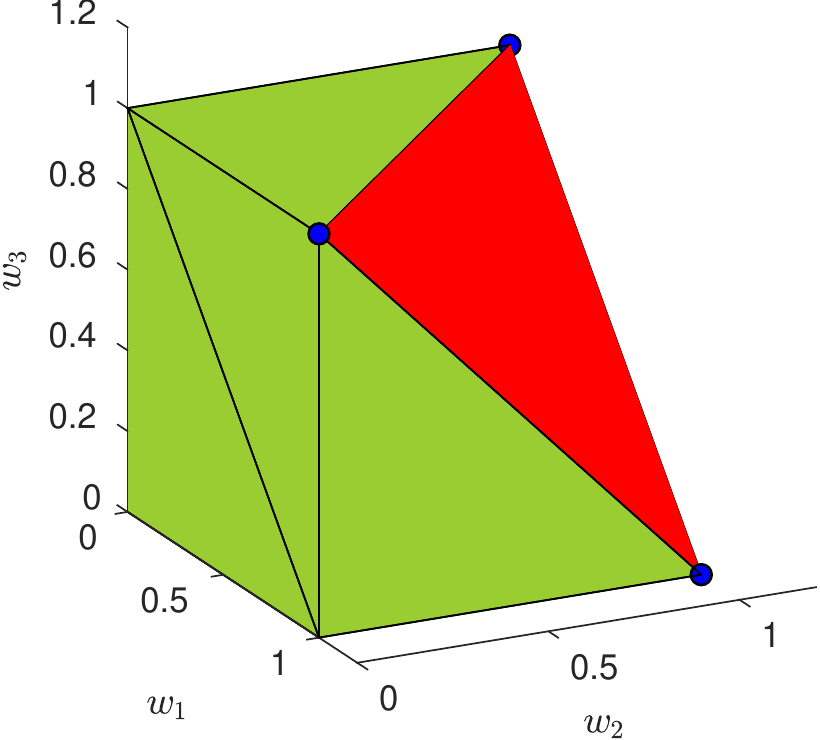}
	\caption{Feasible set \cref{eq:relaxed_OED_problem} for weights $(w_1,w_2,w_3) \in \R^3$ (shown in green) and the restricted simplex~$\Delta_{2}$ (shown in red).}
\end{figure}

Notice that \eqref{eq:relaxed_OED_problem} is a convex problem since the A-criterion is convex w.r.t.\ $\Upsilon$ and $\Upsilon$ depends linearly on $w$; see for instance \cite[Chap.~2.3]{FedorovLeonov2014:1}, \cite[Ch.~10.4.2]{Seber2008}.
In the absence of an upper bound $C_w$ on the total weight, all $w_{k,\ell}$ would be equal to one in an optimal solution since the activation of any of the sensors at any time adds information.
Due to the continuity of the objective and the compactness of the feasible set, it is evident that an optimal distribution of weights exists.
The minimal objective value will be finite since any non-zero weight vector induces a positive definite FIM.
This can be shown using the theory of parabolic equations.
It is easy to see that, without loss of generality, the constraint $\sum_{k,\ell} w_{k,\ell} \le C_w$ will be active at an optimal solution.
A typical picture of the feasible set in three dimensions is shown in \cref{fig:restricted_simplex}, and the red facet represents the subset $\sum_{k,\ell} w_{k,\ell} = C_w$ on which we seek a solution. 
While \cref{fig:restricted_simplex} is only three-dimensional, the (restricted) simplex
\begin{equation}
	\label{eq:simplex}
	\Delta_{C_w}
	\coloneqq
	\Bigg\{ w \in \R^{\Nobs \cdot \Ntime} : 0 \le w_{k,\ell} \le 1, \; \sum_{k,\ell} w_{k,\ell} = C_w \Bigg\} 
\end{equation}
will be of dimension $\Nobs \cdot \Ntime - 1$ and thus rather high dimensional in our numerical experiments.
We thus utilize a simplicial decomposition approach, which restricts problem \eqref{eq:relaxed_OED_problem} to the convex hull of a selection of active vertices of \eqref{eq:simplex}, i.e., a subset of the vertices shown in blue in \cref{fig:restricted_simplex}. 
Simplicial decomposition then takes turns selecting the active vertices and solving problem \eqref{eq:relaxed_OED_problem} on the lower-dimensional simplex spanned by these vertices, which is a subset of $\Delta_{C_w}$.
We refer the reader to \cite[Chap.~9]{Patriksson1999,Patriksson2009} and \cite[Chap.~4]{Bertsekas2015} for a general account on simplicial decomposition.
By \cite[Prop.~4.2.1]{Bertsekas2015}, a solution to \eqref{eq:relaxed_OED_problem} will be found in finitely many iterations.

In each iteration, one active vertex is added based on the magnitude of the partial derivative of the objective. 
Subsequently, problem \eqref{eq:relaxed_OED_problem} --- restricted to the updated lower dimensional simplex --- is solved via Torsney's algorithm, which is a simple iterative scheme to update the weights.
We refer the reader, e.g., to \cite[Ch.~4]{HerzogRiedelUcinski2017:1} for details.
Both steps of the algorithm make use of the following result concerning the derivative of the objective.

\begin{lemma}
	\label{lemma:derivative_of_OED_objective}
	Suppose that $w$ is a weight vector such that $\Upsilon(w)$ is positive definite.
	Then the objective in \eqref{eq:relaxed_OED_problem} is differentiable w.r.t.\ $w$.
	Its partial derivative w.r.t.\ $w_{k,\ell}$ is given by 
	\begin{equation}
		\label{eq:partial_derivative_of_OED_objective}
		\frac{\partial}{\partial w_{k,\ell}} \Phi_A(\Upsilon(w))
		=
		- \trace \left( \Upsilon(w)^{-1} B \, \Upsilon(w)^{-1} \Upsilon_{k,\ell} \right)
		.
	\end{equation}
\end{lemma}
This result can be found, for instance, in \cite[Thm.~B.19]{Ucinski2005:1} and \cite[Ch.~17.8]{Seber2008}. 

In our implementation, we pre-calculate the elementary FIMs $\Upsilon_{k,\ell}$ \eqref{eq:elementary_FIM}.
This requires the solution of one forward problem \eqref{eq:forward_problem} as well as a number of sensitivity equations \eqref{eq:sensitivity_problem_weak}, depending on the dimension~$\Nbasis$ of the subspace of parameter variations.
Recall that the full FIM associated with a weight vector $w$ is given by \eqref{eq:combined_experiment}.

A necessary and sufficient optimality condition for \eqref{eq:relaxed_OED_problem} is that
\begin{equation}
	\label{eq:optimality_condition_for_relaxed_OED_problem}
	- \frac{\partial}{\partial w_{k,\ell}} \Phi_A(\Upsilon(w))
	\begin{cases}
		\ge \xi & \text{if } w_{k,\ell} = 1, \\
		= \xi & \text{if } 0 < w_{k,\ell} < 1, \\
		\le \xi & \text{if } w_{k,\ell} = 0.
	\end{cases}
\end{equation}
holds for some positive $\xi$; see \cite[Prop.~1]{UcinskiPatan2007} for a similar OED criterion and \cite[Thm.~1]{Pronzato2004} for a general result.
An inexact version of \eqref{eq:optimality_condition_for_relaxed_OED_problem} will be used as the stopping criterion for the simplicial decomposition iteration in our numerical experiments. 

\section{Numerical Results}
\label{sec:Numerical_results}

In this section we present and discuss numerical results showing optimal sensor activation patterns for problem \eqref{eq:relaxed_OED_problem}.
We concentrate on three examples, one in 2D and two in 3D, which differ further by the choice of boundary conditions and sensor types.
In all situations, the underlying model is given by the parabolic diffusion equation~\eqref{eq:forward_problem}.
The weak formulation of both the forward problem~\eqref{eq:forward_problem_weak} and the sensitivity equation~\eqref{eq:sensitivity_problem_weak} are discretized using linear finite elements on triangles and tetrahedra, respectively.
Time is discretized using the backward Euler time stepping on a uniform time grid.
We use the finite element toolbox GetFEM++ \cite{GETFEM}.
The diffusion coefficient is chosen to be $\kinc = \num{1e-3}$ inside the inclusion $\Omegainc$ and $\kbulk = \num{1e-1}$ in $\Omega\setminus \mathrm{cl}(\Omegainc)$.
In all numerical examples we have Dirichlet data $u_D = 1$ and initial conditions $u(\cdot, 0) = 0$.
The domains are chosen to be the unit cubes, i.e., $\Omega = (0,1)^d$ for $d=2,3$.
All meshes are generated with GMSH \cite{GeuzaineRemacle2009}.
In both two and three dimensions, we extend normal vector fields $\bar V$ representing shape variations to deformation fields $V$ in the hold-all $D$ by solving \eqref{eq:elasticity_problem} with the Lamé parameter $\lambda = 0.01$ and $\mu = 0.495$.

\subsection{The 2D Case}
\label{subsec:Numerical_results_2D}

We first discuss the two dimensional example, where the measurement regions $\Omegaobs^k$ are assumed to be two dimensional subregions of $\Omega = (0,1)^2$.
Here the inclusion $\Omegainc$ is modeled by a B-spline and it is located inside the hold-all $D = (0.35, 0.65)^2$.
$\Nobs = 8$ sensors of equal size are positioned around $D$, each of them given by a square with edge length $0.3$; see \cref{fig:domain_experiment}.
The time horizon of the experiment is chosen as $T = 10$, which is discretized into 21~equal intervals.
We assume that all time steps are potential measurements, which leads to $\Ntime = 22$ in the optimum experimental design algorithm.
Consequently, we have $\Nobs \cdot \Ntime = 176$~potential sensor activations.
The upper bound for the sum of the weights in the optimum experimental design is chosen as $C_w = 10$.
In equation \eqref{eq:covariance_root_differential_operator} the parameters which encode the variation and correlation length within one sensor region are chosen as $\alpha_0 = 0.01$ and $\alpha_1=1.0$.

The upper segment of the boundary serves as Dirichlet boundary $\GammaD$, and Robin boundary conditions \eqref{eq:forward_problem} are imposed on the remainder $\GammaR$.
We choose $\beta = 10.0$ in the left half of the bottom part of the boundary, and $\beta = 0$ elsewhere on $\GammaR$.

The shape deformations $V$ are taken from the nine dimensional space spanned by the basis functions $\bar V_i$ depicted in \cref{fig:untransformed_bumps}, which are then extended to functions $V_i$ on $D$ by solving \eqref{eq:elasticity_problem}.
In our numerical experiments we choose one of the basis functions equal to a constant normal vector field along the interface $\Gammainc$; see the last plot in \cref{fig:untransformed_bumps}.
The other basis functions $\bar V_i$ are chosen to be equidistantly distributed Gaussians on $\Gammainc$.
To achieve this we proceed as follows.
Let $c: [0,L) \to \R$ be a curve representing $\Gammainc$ and $L$ its length.
We choose a parametrization by arc length and can thereby define nodes on $\Gammainc$ equidistantly by their parametrization $r_i \in [0, L)$.
Since we assume an arbitrarily discretized curve representing $\Gammainc$, we recover $L$ be adding the length of all discrete boundary segments.
The parametrization is then found by marking one surface node as starting point and iteratively adding the lengths of connected segments.
This enables us to choose the boundary deformation as
\begin{equation}
	\bar V_i(r) = n(r) \cdot \exp(-s\cdot d(r, r_i)^2)
	\label{eq:gaussian_bump}
\end{equation}
for all $r \in [0, L)$ where $n(r)$ denotes the unit outer normal vector at position $r$ and $d_i(r, r_i) = \min\big(\vert r - r_i\vert, \min(L - r + r_i, L-r_i +r)\big)$ approximates the geodesic distance.
Further, $s$ denotes a slope factor which controls how fast $\bar V_i$ decreases and it is chosen as $s = 100$ for the 2D example.

\begin{figure}[htbp]
	\label{fig:final_weights_and_diffusion}
	\begin{center}
		\setlength{\fboxsep}{4pt}%
		\setlength{\fboxrule}{1pt}%
		\fbox{\includegraphics[width=.6\textwidth]{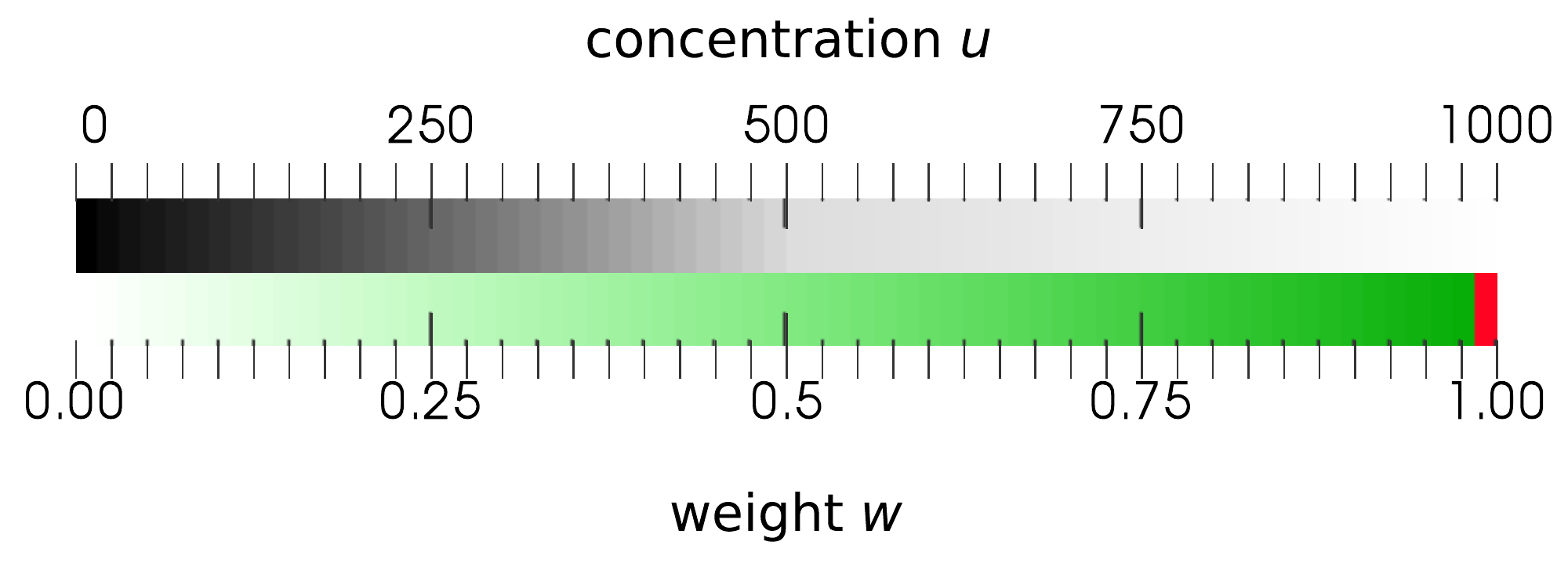}}\\
		
		\vspace{0.5em}
		\newcolumntype{C}{>{\centering\arraybackslash} m{0.3\textwidth} }
		\begin{tabular}{CCC}
			\toprule
			\includegraphics[width=.3\textwidth]{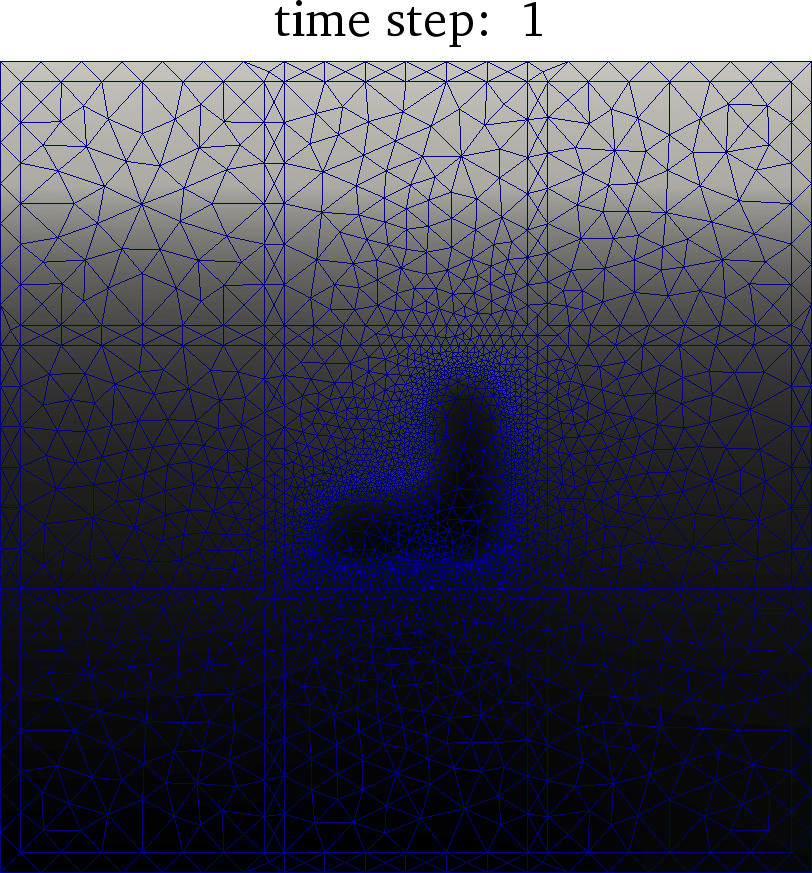}&
			\includegraphics[width=.3\textwidth]{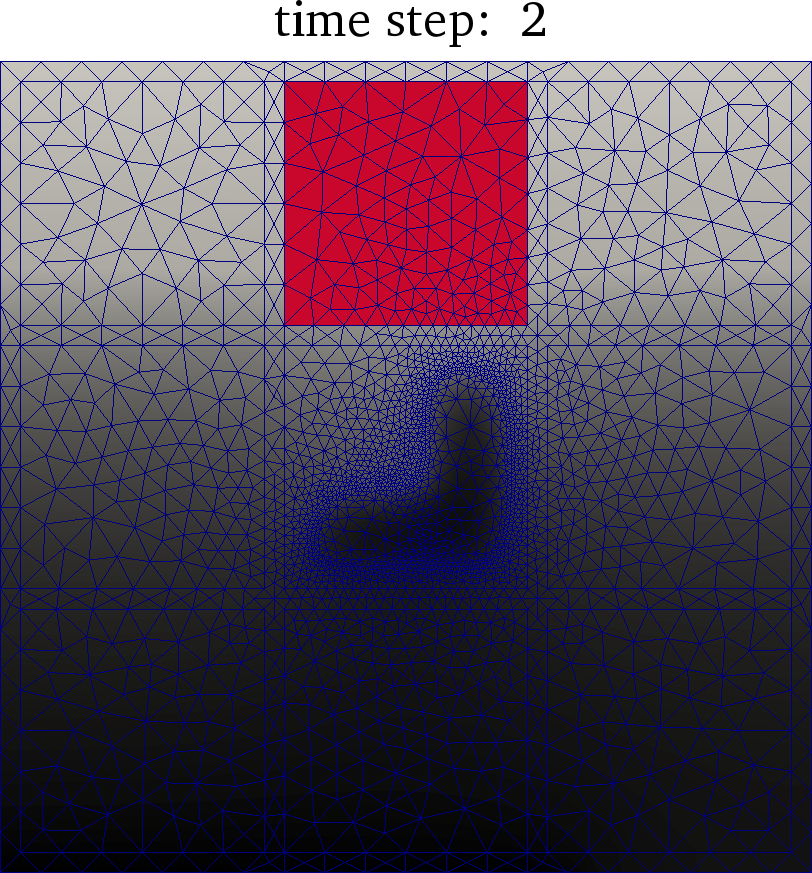}&
			\includegraphics[width=.3\textwidth]{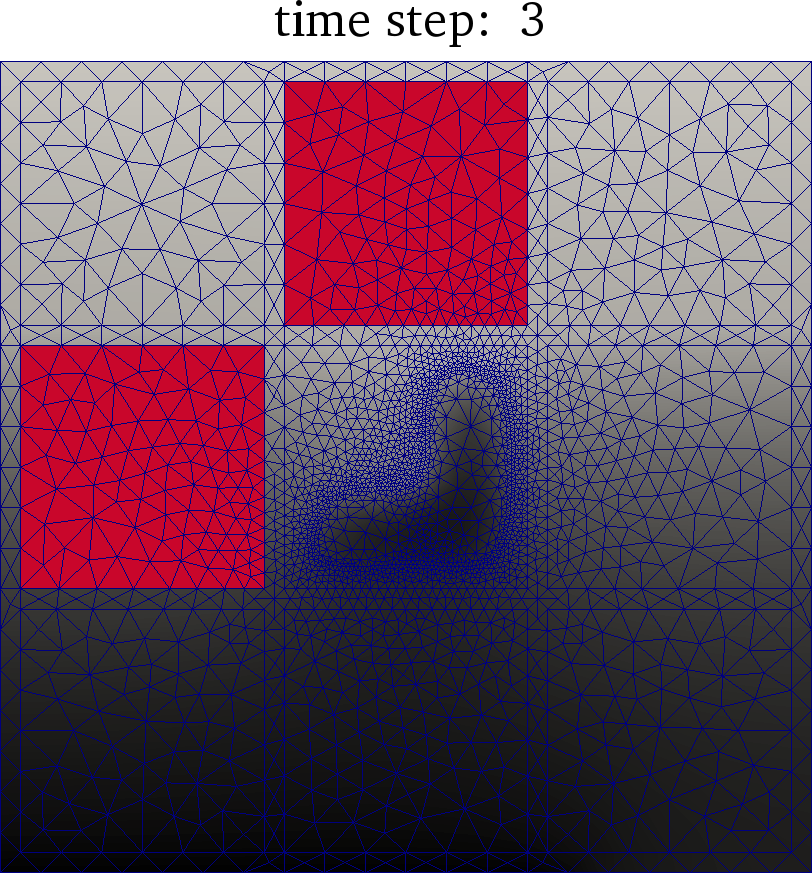}\\
			\midrule
			\includegraphics[width=.3\textwidth]{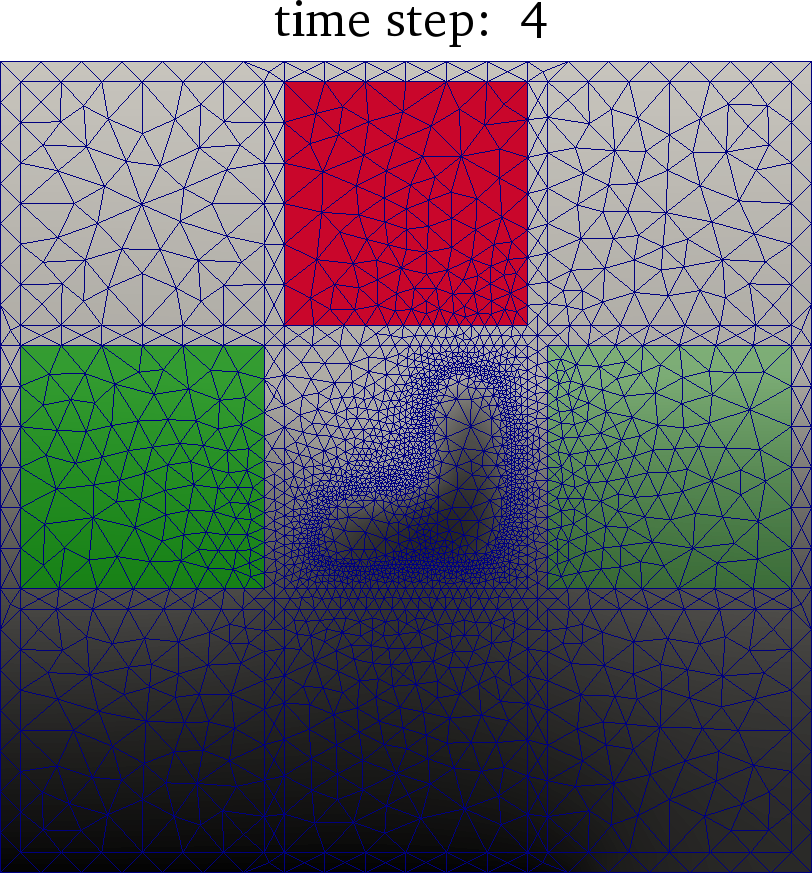}&
			\includegraphics[width=.3\textwidth]{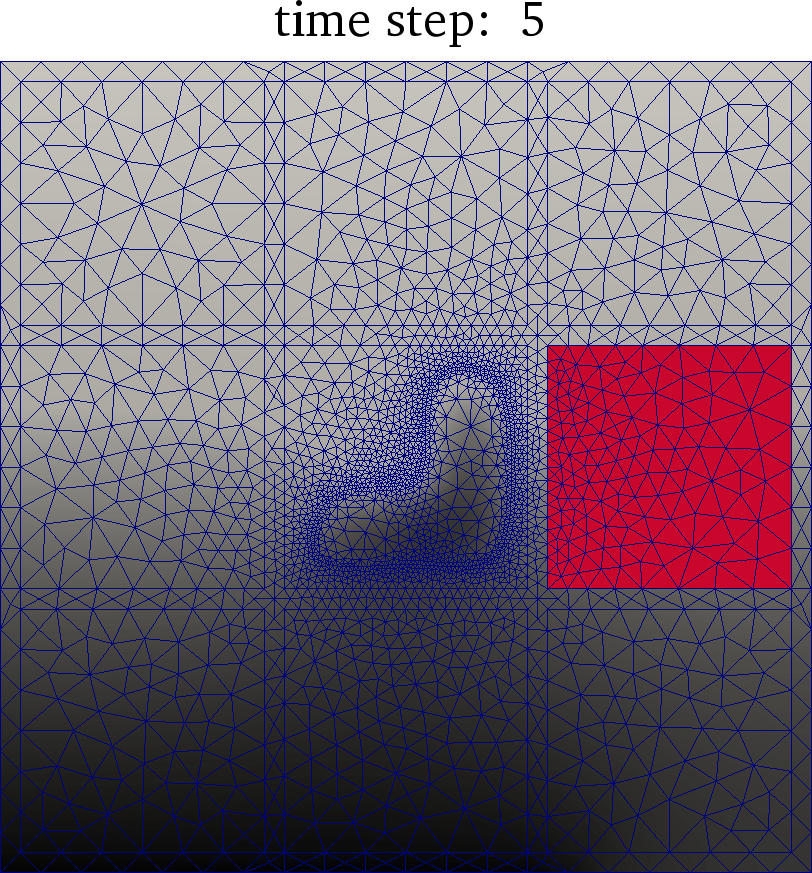}&
			\begin{minipage}{.3\textwidth}
				no activation from time step~6 until 14
			\end{minipage} 
			\\
			\midrule
			\includegraphics[width=.3\textwidth]{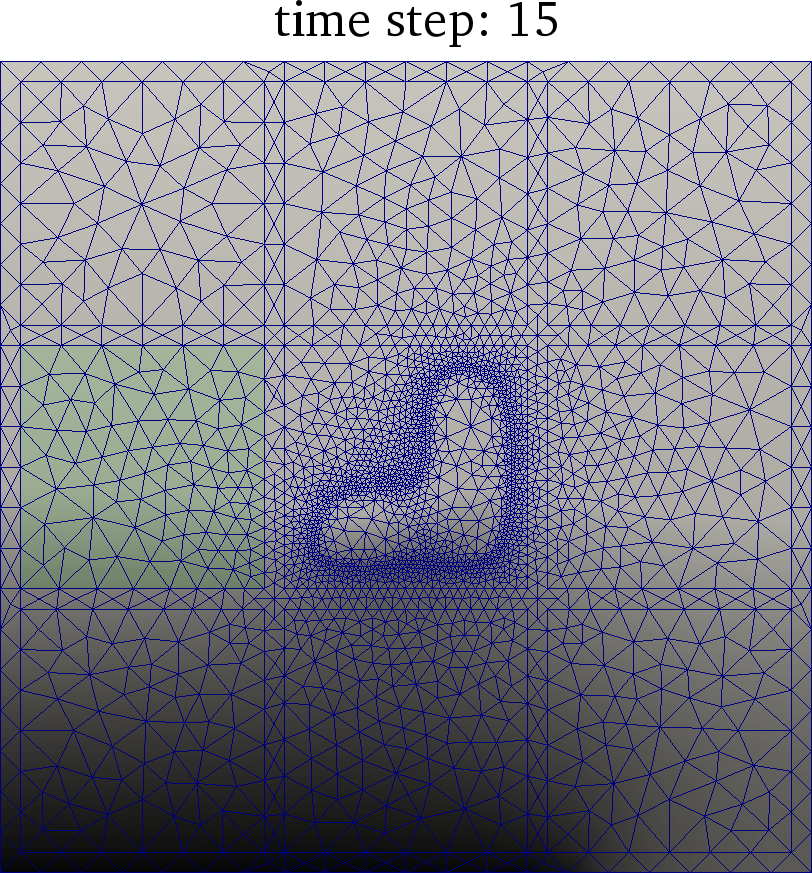}&
			\includegraphics[width=.3\textwidth]{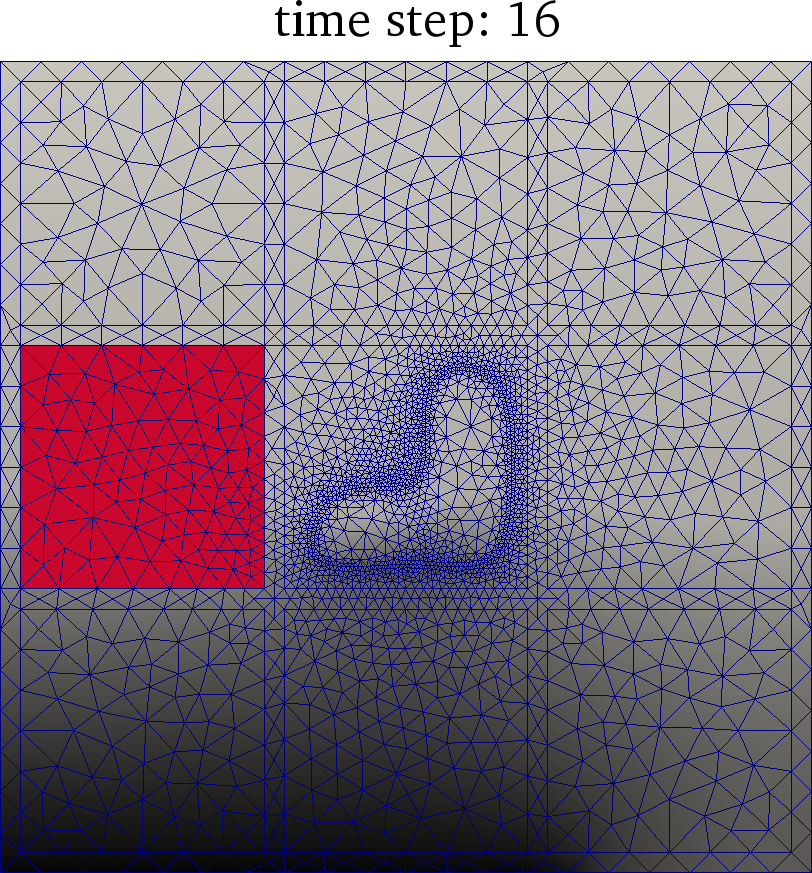}&
			\begin{minipage}{.3\textwidth}
				no activation from time step~17 until 18
			\end{minipage} 
			\\
			\midrule
			\includegraphics[width=.3\textwidth]{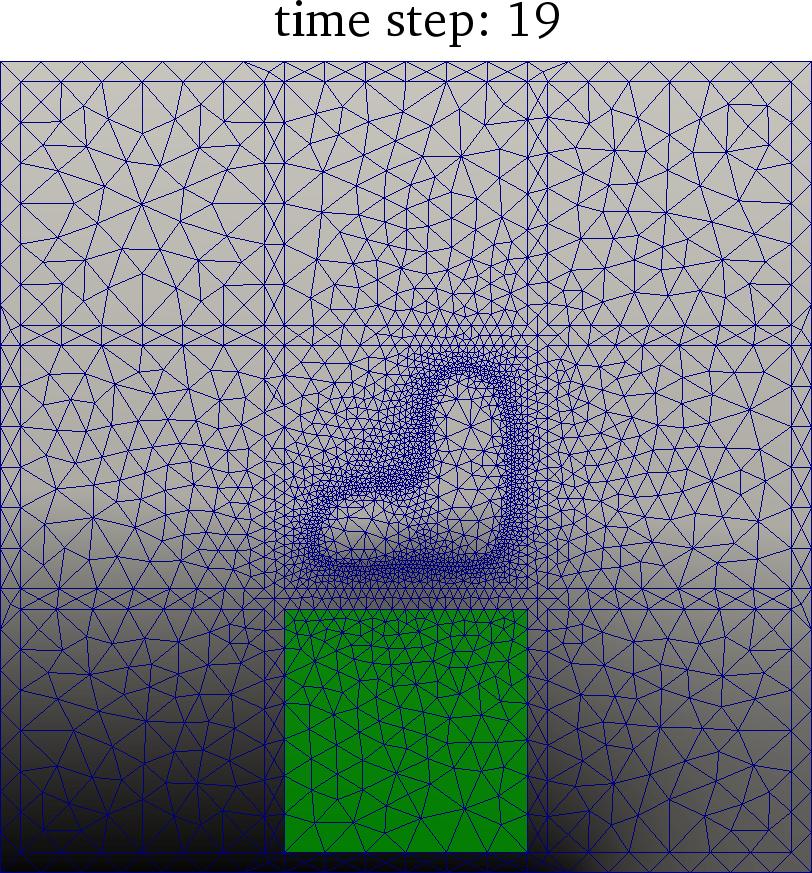}&
			\includegraphics[width=.3\textwidth]{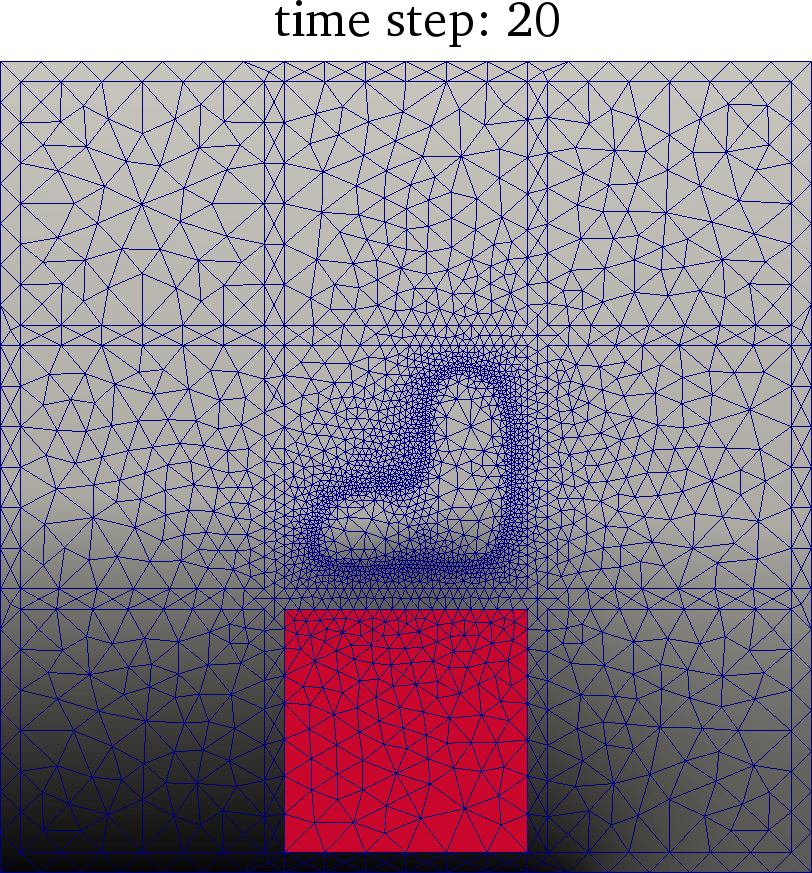}&
			\includegraphics[width=.3\textwidth]{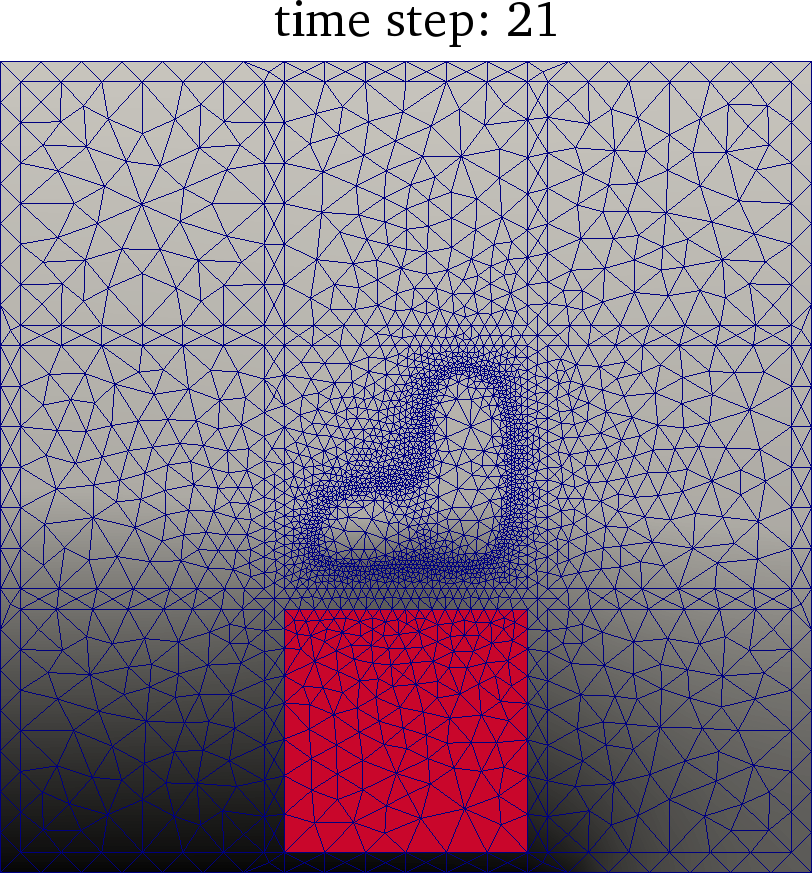}\\
			\bottomrule
		\end{tabular}
		\caption{Forward solution showing the concentration in gray-scale according to \eqref{eq:forward_problem}, and optimal sensor activation pattern in red--green color scale.}
	\end{center}
\end{figure}

The optimal activation for the sensor pattern as an overlay over the diffusion process is shown in \cref{fig:final_weights_and_diffusion}. 
Notice that sensors tend to get activated mainly when the diffusion front arrives, since this event gives rise to significant spatial gradients.
This is reflected in the sensitivies (not shown) according to \eqref{eq:sensitivity_problem_weak}.

\begin{figure}[htbp]
	\begin{center}
		\begin{tabular}{cc}
			\includegraphics[width=.45\textwidth]{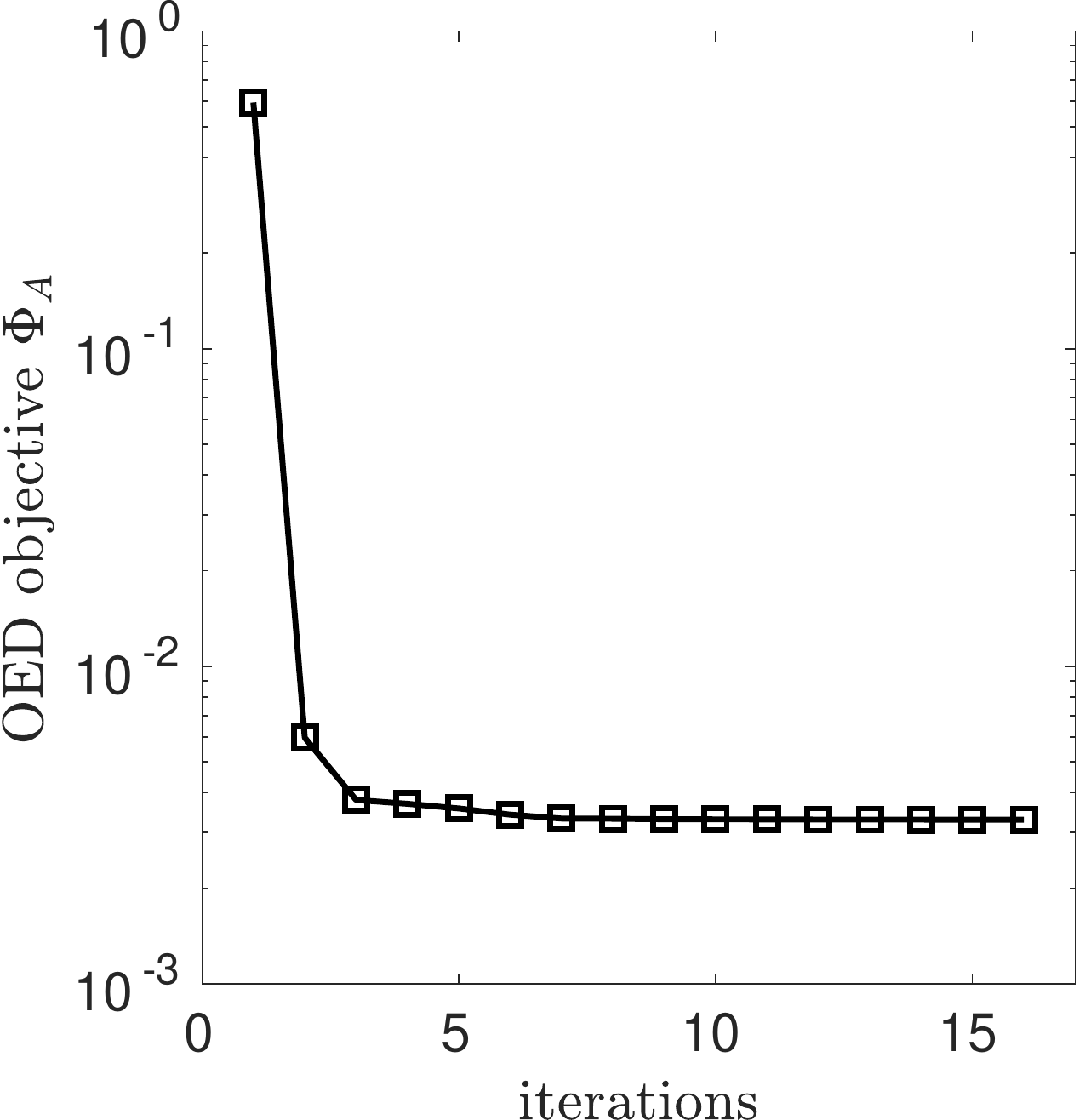}&
			\includegraphics[width=.45\textwidth]{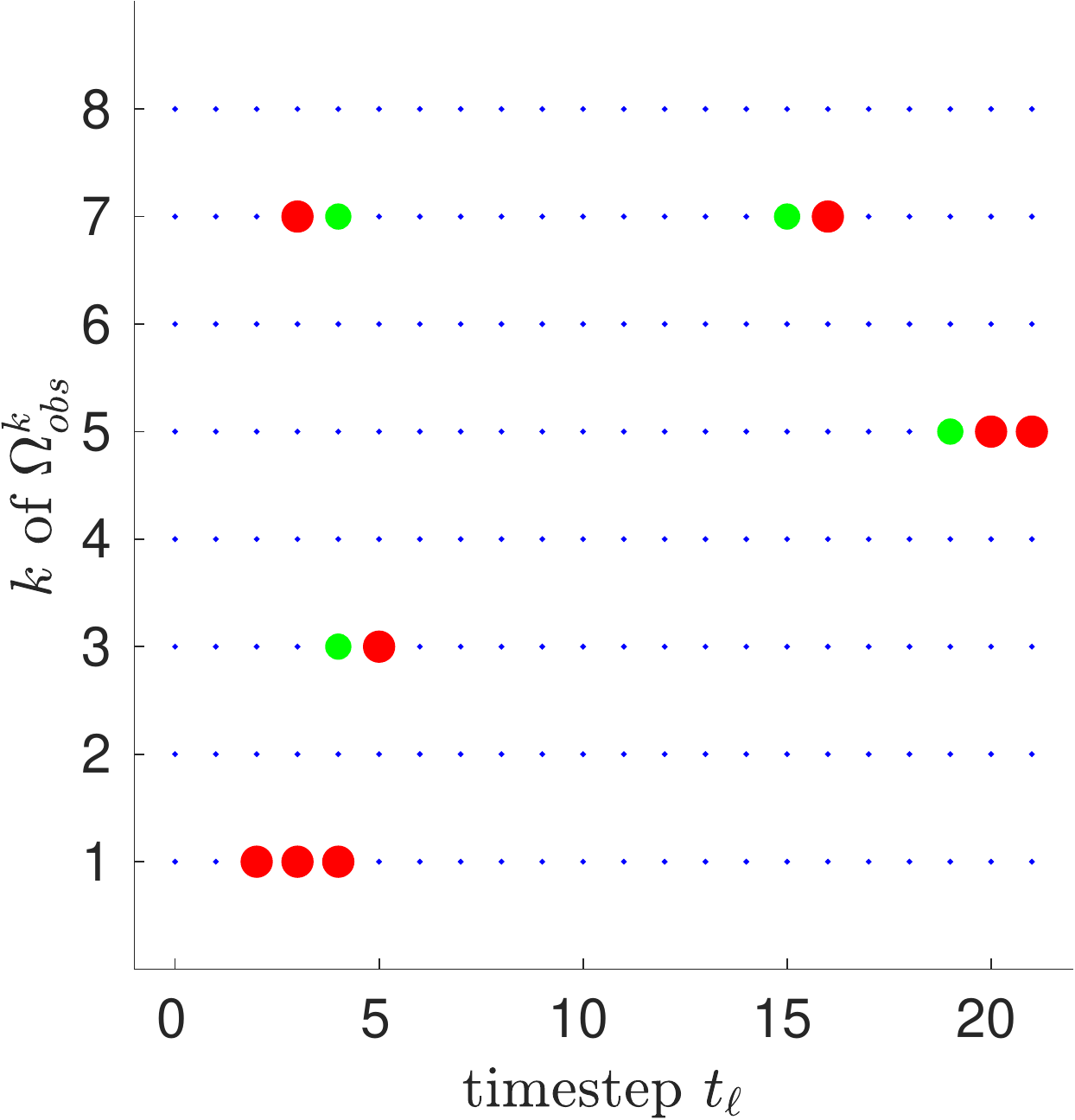}\\
			\includegraphics[width=.45\textwidth]{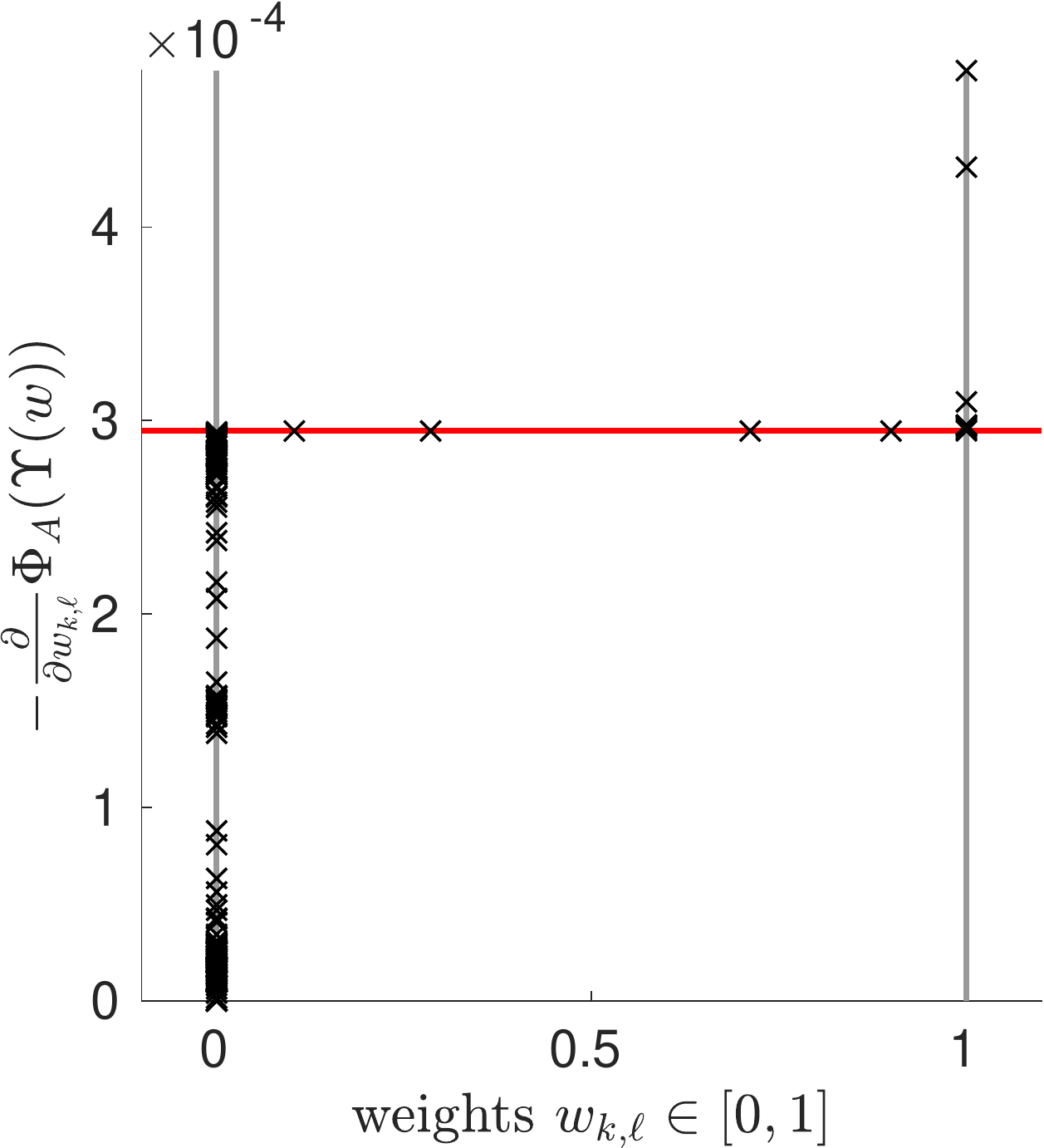}&
			\includegraphics[width=.45\textwidth]{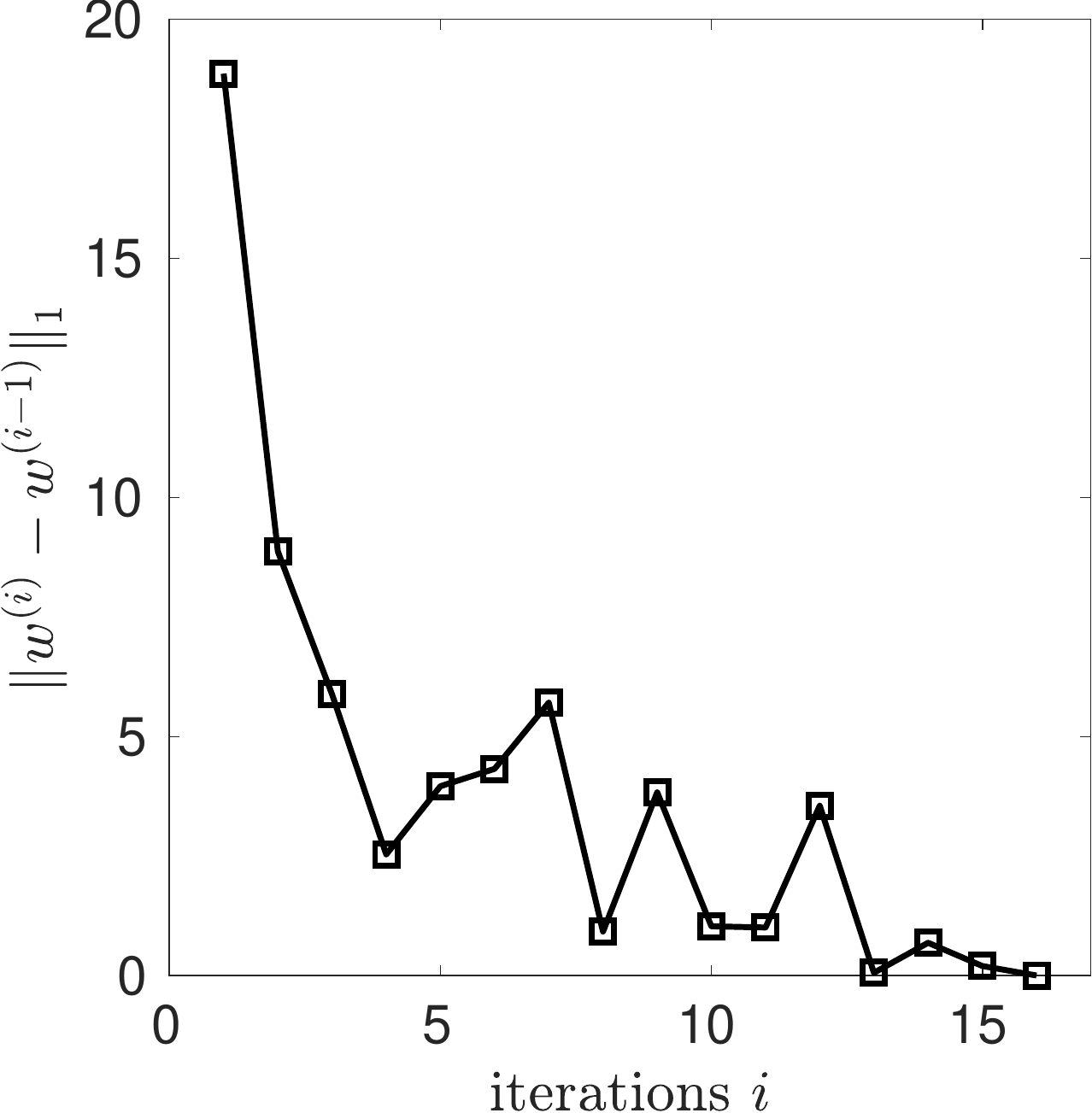}
		\end{tabular}
	\end{center}
	\caption{Results for the 2D experiment: evolution of the objective $\Phi_A$ in a semi-log scale (top left); optimal sensor activation pattern with $w_{k,\ell}=1$ shown in red and $w_{k,\ell} \in \left(0,1\right)$ shown in green. The measurement domains $\Omegaobs^k$ are numbered as in \cref{fig:domain_experiment} (top right); the optimality condition \eqref{eq:optimality_condition_for_relaxed_OED_problem} is verified at the final iterate (bottom left); evolution of the change for the weights measured in the 1-norm (bottom right).}
	\label{fig:results_2D}
\end{figure}

Owing to the simplicial decomposition approach we obtain a rather sparse pattern with $8$~weights equal to one, $4$~weights in $\left(0,1\right)$ and all the other $164$~weights exactly equal to zero. 
The sparsity of the activation pattern can also be seen in the top right picture of \cref{fig:results_2D}, where $w_{k,\ell}=1$ is indicated by red dots and $w_{k,\ell}\in\left(0,1\right)$ is shown in green.
The iteration history is depicted in \cref{fig:results_2D}, where we started from an initial guess of uniform activation of all weights throughout space and time.
The necessary and sufficient optimality condition \eqref{eq:optimality_condition_for_relaxed_OED_problem} is fulfilled after $16$~iterations of the outer loop in the simplicial decomposition algorithm, see the bottom left picture in \cref{fig:results_2D}, where the red horizontal line marks the value of~$\xi$ in \eqref{eq:optimality_condition_for_relaxed_OED_problem}.
A rapid decrease of the objective function $\Phi_A$ \eqref{eq:relaxed_OED_problem} in the first iterations can be seen in the top left of \cref{fig:results_2D}.
The subsequent iterations mainly serve to improve the sparsity pattern of the sensor activation.
The change of the weight vector between two subsequent iterations, measured in the 1-norm, is shown in the bottom right of \cref{fig:results_2D}.
Also here, the largest changes appear in the first few iterations. 
Due to the relaxation of the weights $w_{k,\ell}$ in \eqref{eq:relaxed_OED_problem} we cannot expect to obtain a truly binary solution.
If desired, a rounding heuristic can be applied, such as setting the largest $C_w$ weights to one and the remaining ones to zero.

\begin{figure}[htbp]
	\label{fig:optimized_bumps_2d}
	\begin{center}
		\begin{tabular}{ccc}
			\includegraphics[width=.3\textwidth]{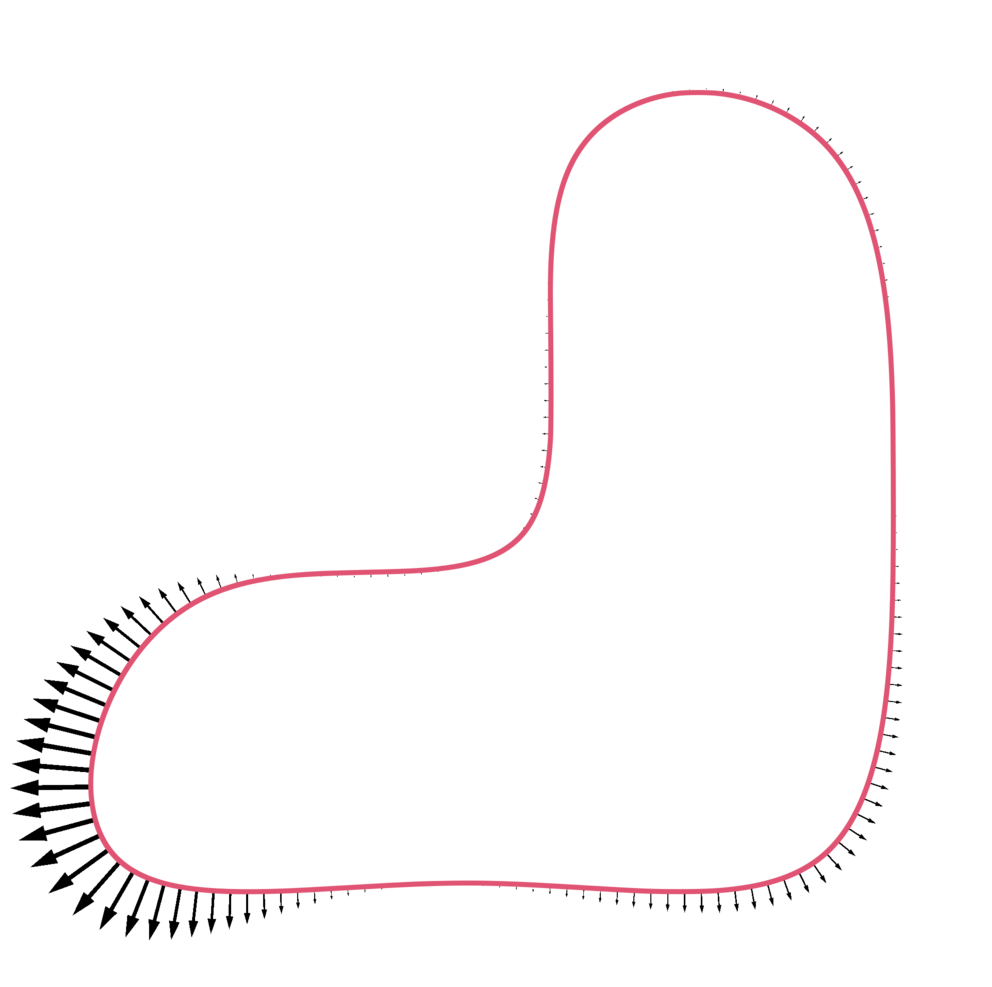}&
			\includegraphics[width=.3\textwidth]{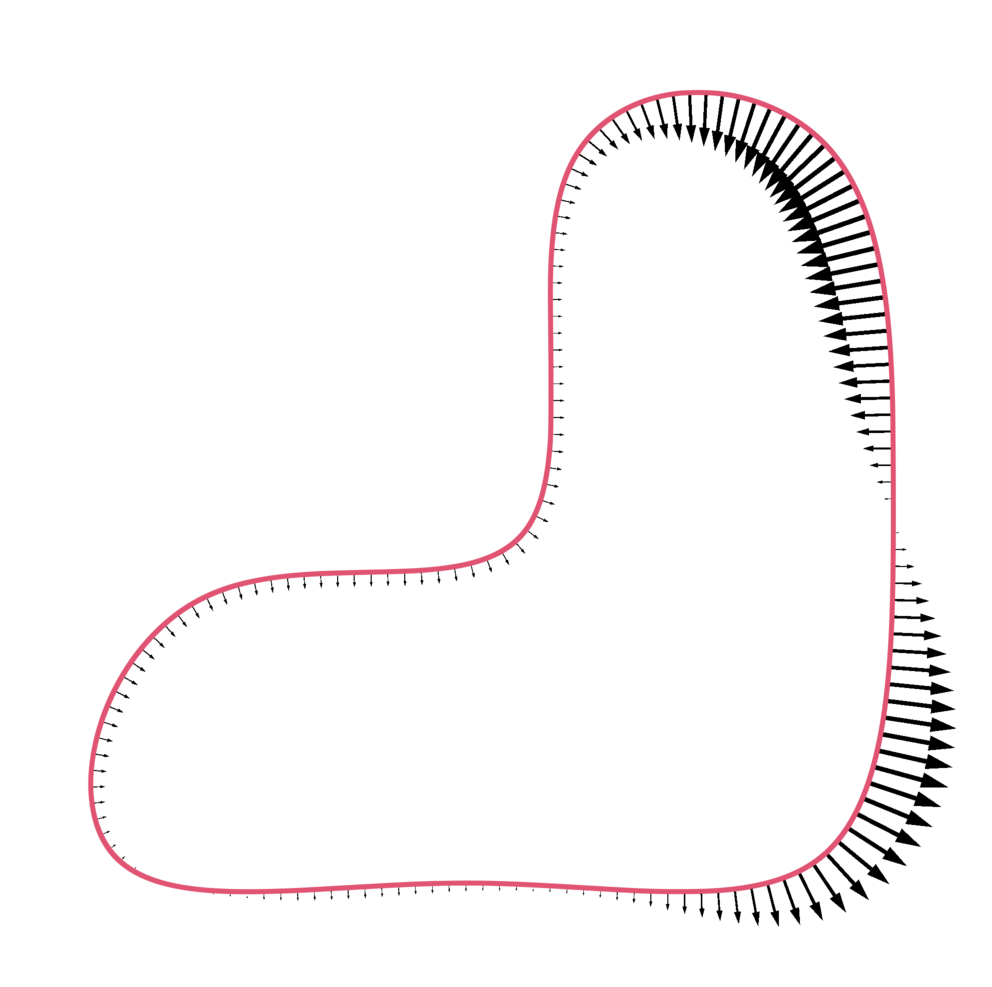}&
			\includegraphics[width=.3\textwidth]{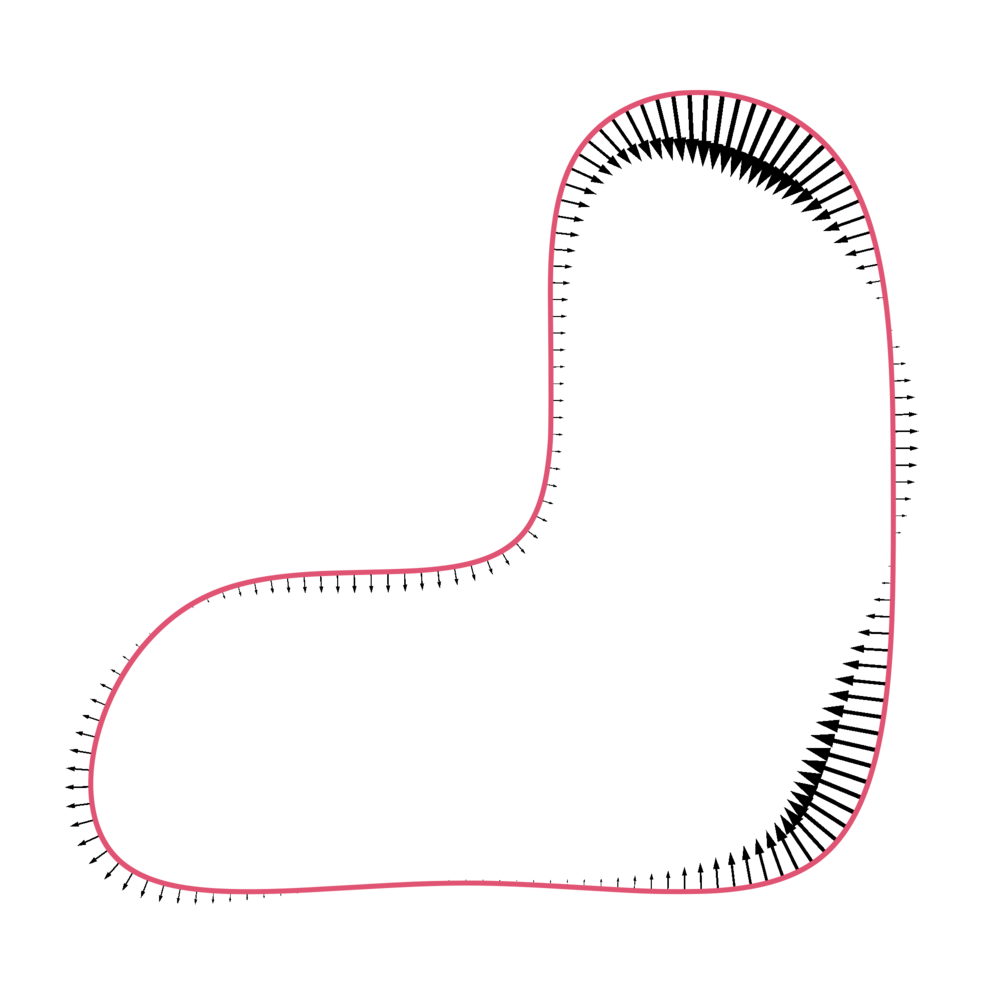}\\
			\includegraphics[width=.3\textwidth]{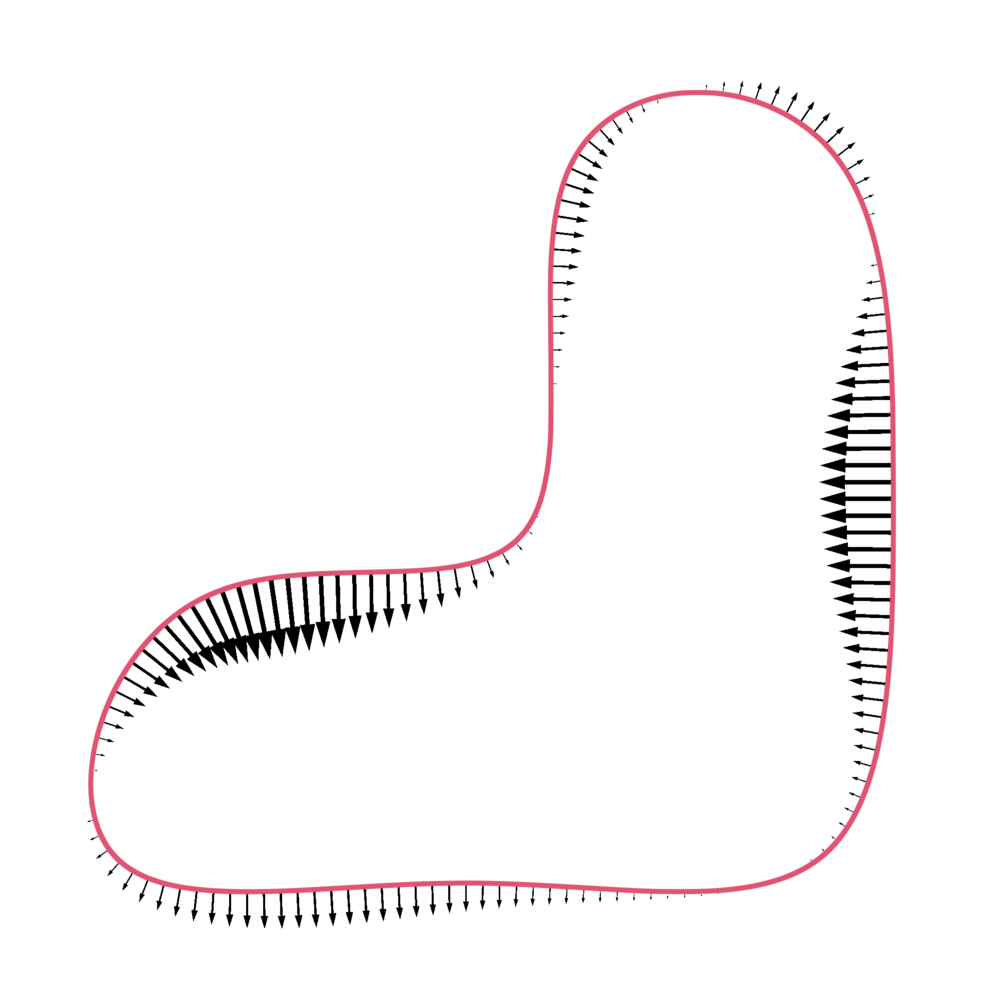}&
			\includegraphics[width=.3\textwidth]{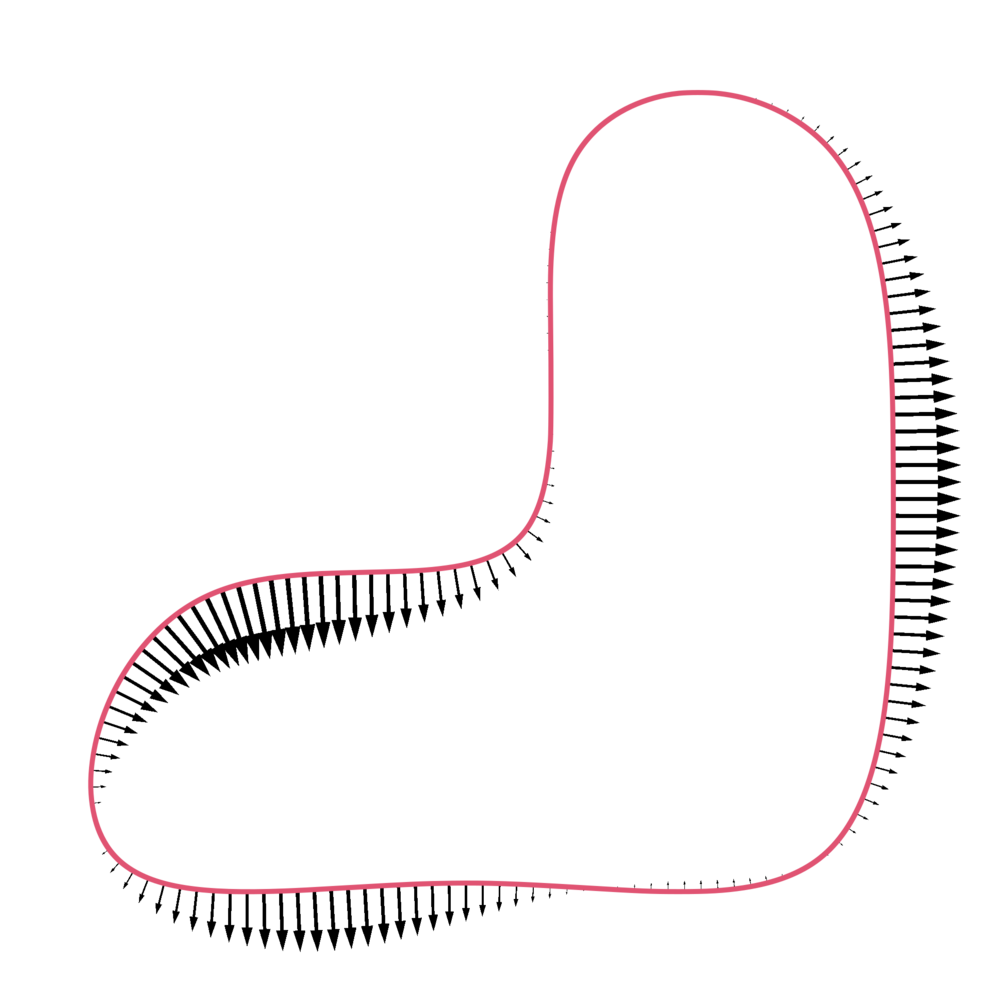}&
			\includegraphics[width=.3\textwidth]{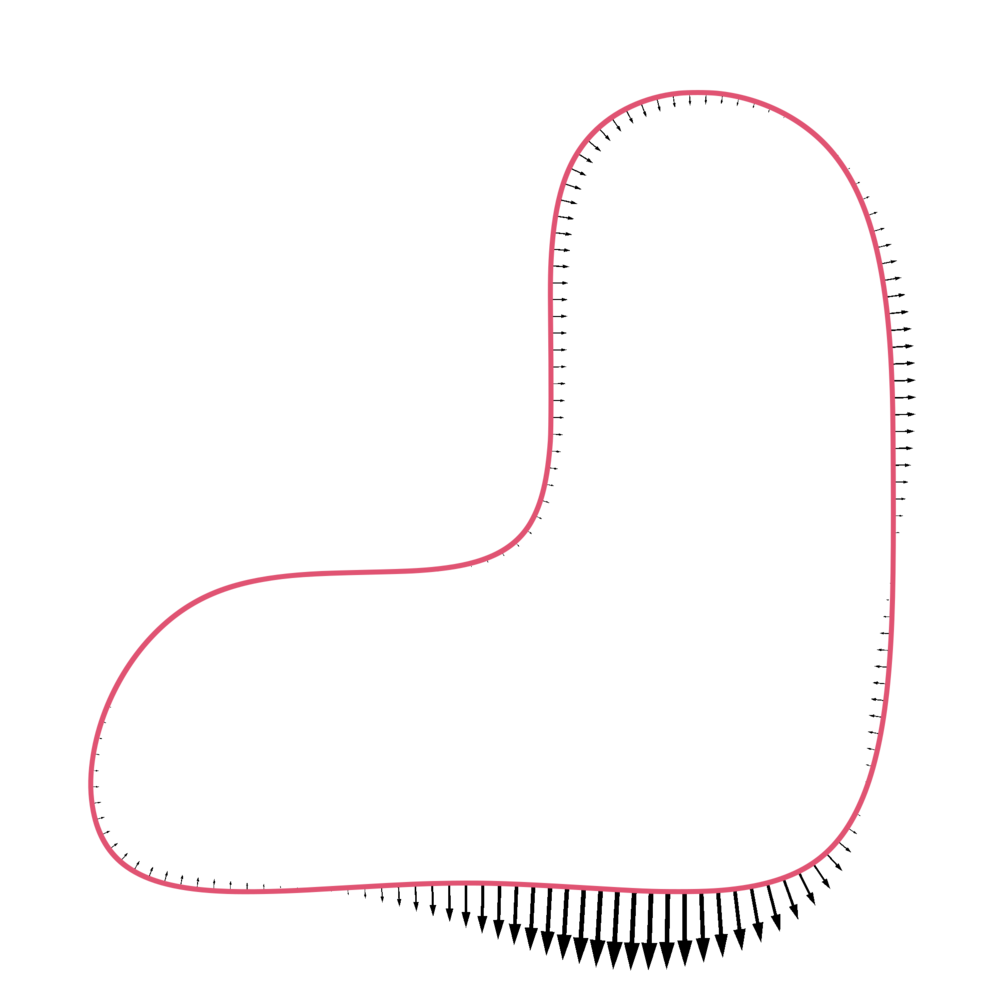}\\
			\includegraphics[width=.3\textwidth]{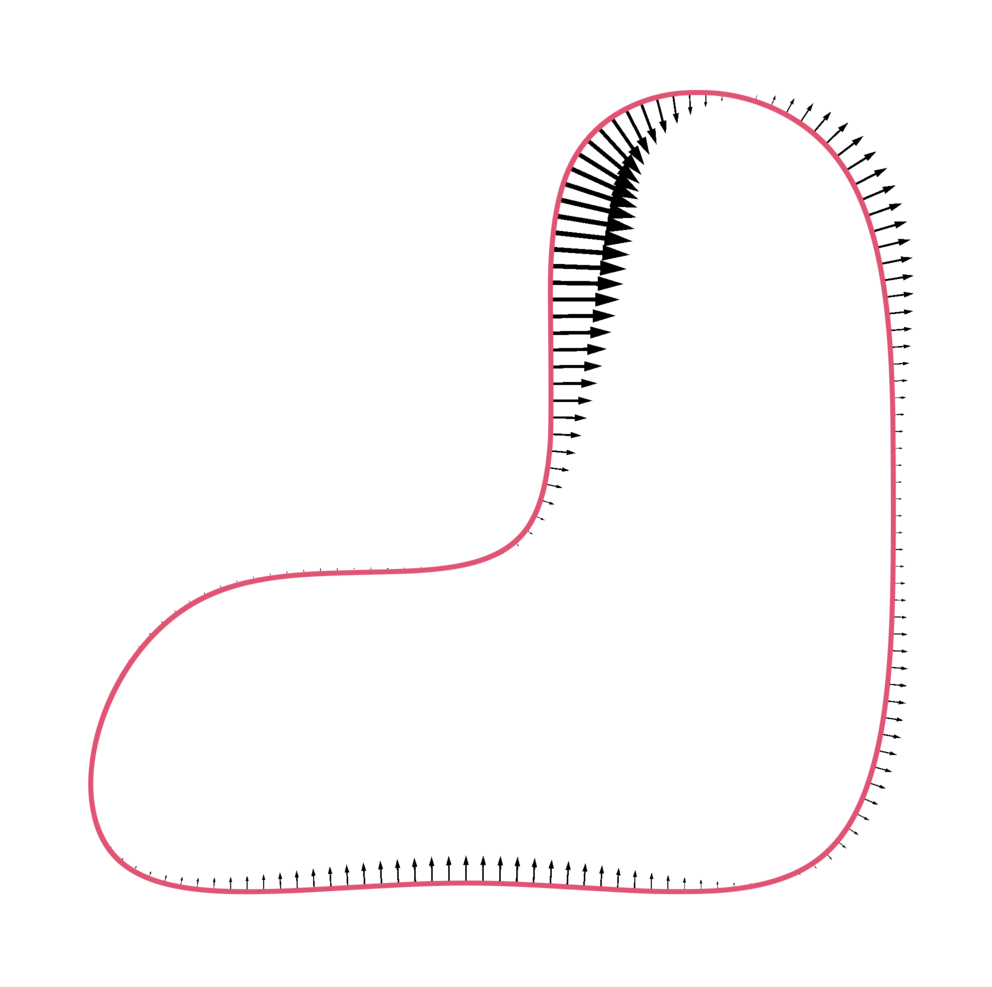}&
			\includegraphics[width=.3\textwidth]{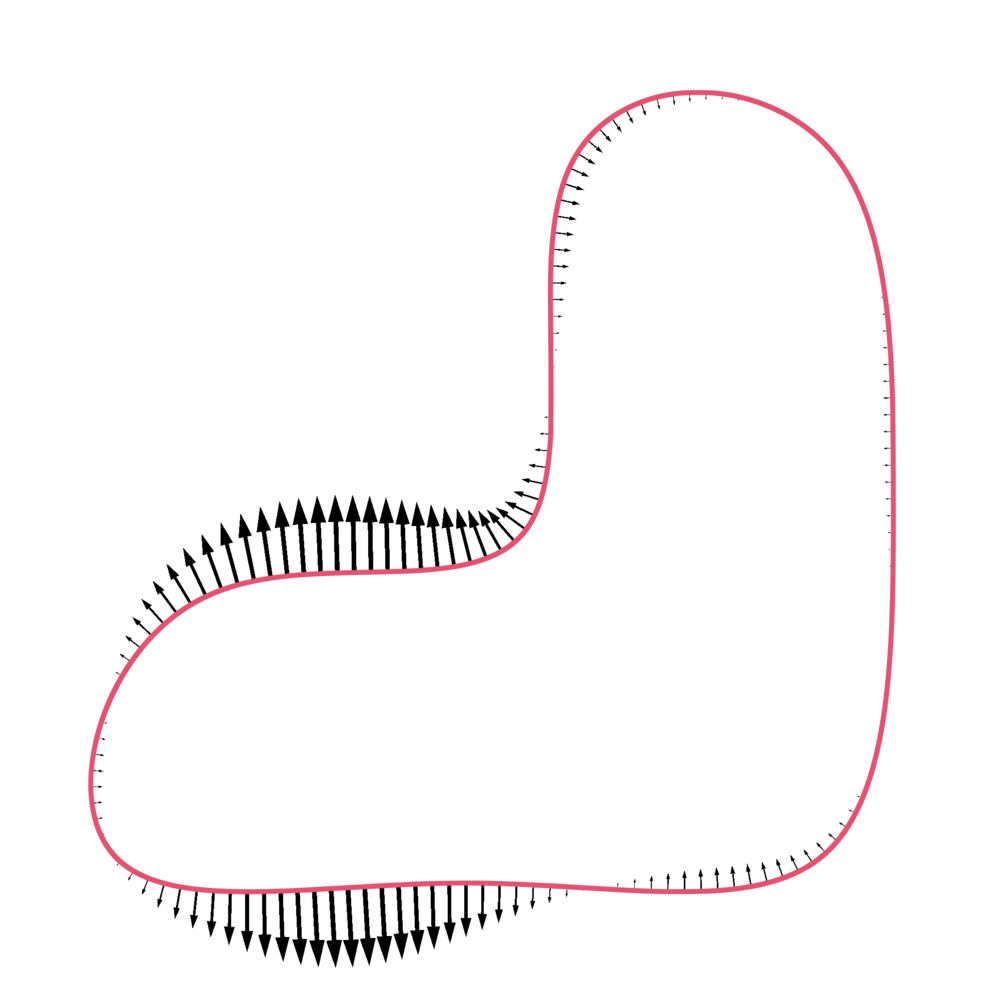}&
			\includegraphics[width=.3\textwidth]{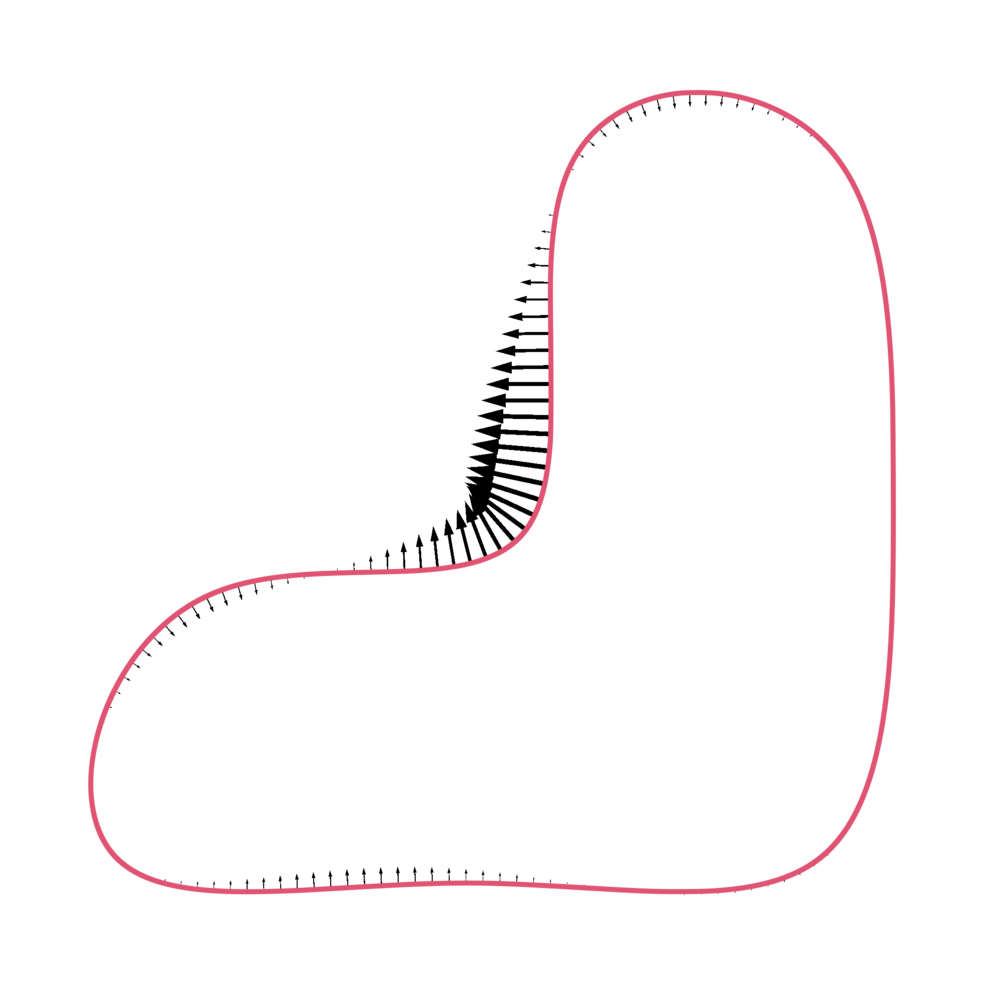}
		\end{tabular}
		\caption{Eigenfunctions (orthonormal with respect to $B$ from \eqref{eq:generalized_eigenvalue_problem}, but shown here with unified lengths) according to the optimized weights. The corresponding eigenvalues $(\Lambda_1^{-1},\ldots,\Lambda_9^{-1})$ are given from top left to bottom right by 
		$(0.27,
		0.49,
		0.77,
		1.04,
		1.24,
		2.58,
		3.39,
		6.68,
		16.48)\times\num{e-4}$. 
		}
	\end{center}
\end{figure}

Once an optimal weight vector $w^\ast$ and the associated optimal FIM $\Upsilon(w^\ast)$ are determined, we solve the generalized eigenvalue problem \eqref{eq:generalized_eigenvalue_problem}.
The eigenfunctions are shown in \cref{fig:optimized_bumps_2d} and the reciprocals of the eigenvalues are listed in the caption of that figure.
Recall that the sum of the reciprocal eigenvalues equals the value of the objective function $\Phi_A$; see \eqref{eq:A-criterion}.
The eigenfunctions are sorted row-wise from left to right, such that the top left eigenfunction spans the subspace in the shape deformation tangent space with the best identifiability.
We can interpret \cref{fig:optimized_bumps_2d} in the sense that the elongation of the inclusion $\Omegainc$ to the left can be identified with highest reliability.
That is, variations of $\Omegainc$ in this direction have the largest impact among all perturbations of equal $b$-norm \eqref{eq:Riemannian_metric} on the selected state measurements.
Notice that this is due to the choice of the Robin boundary condition, which allows for a continuous outflow across the bottom left part of the boundary $\GammaR$.

\begin{figure}[htbp]
	\label{fig:optimal_eigenvalues_2d}
	\begin{center}
		{%
			\setlength{\fboxsep}{7pt}%
			\setlength{\fboxrule}{1pt}%
			\fbox{\includegraphics[width=.6\textwidth]{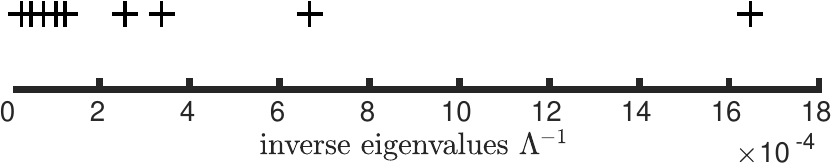}}
		}
		\caption{Inverse eigenvalues $(\Lambda_1^{-1},\ldots,\Lambda_9^{-1})$ according to the FIM of the optimized weights. The eigenvalues are
		$(0.27,
		0.49,
		0.77,
		1.04,
		1.24,
		2.58,
		3.39,
		6.68,
		16.48)\times\num{e-4}$. 
		}
	\end{center}
\end{figure}

By contrast, the eigenfunction in the bottom right visualizes the shape deformation with the worst identifiability, and there is a noticeable gap between the last and next-to-last eigenvalue in \cref{fig:optimal_eigenvalues_2d}.
We attribute this gap at least partially to the fact that the corresponding shape deformation has the largest distance to any sensor position.
It is thus likely that interface perturbations in this direction only have little influence on the measurements.
Note that eigenfunctions of \cref{fig:optimized_bumps_2d} are linear combinations of the basis functions shown in \cref{fig:untransformed_bumps}, and recall that only the subspace spanned by these basis functions is considered.
Notice also that the eigenfunctions, extended to functions on $D$ by \eqref{eq:elasticity_problem}, are orthogonal w.r.t.\ the inner product \eqref{eq:Riemannian_metric}.

In what follows we compare a number of settings described in \cref{tab:compare_2D_settings}, which differ with respect to the transfer coefficient $\beta$ on the boundary $\GammaR$.
We also compare optimized weights with uniform weights.
All other settings are as described above.
The corresponding value of the objective $\Phi_A$ and the inverse eigenvalues of the resulting FIM are displayed in \cref{tab:compare_results_2D_settings}. 
The best possible identification result overall is obtained for case~1, which we discussed in detail above. 
For the second case, where we choose $\beta=10$ at the bottom \emph{right} half of $\GammaR$, we see a similar distribution of the inverse eigenvalues and a slightly increased objective value.
In case~3 we consider the unoptimized version of case~1 with equally distributed measurement weights of the same total mass $C_w$ as before.
The inverse eigenvalues each increase by a factor approximately between~2 and 5.
This demonstrates the gain of information possible through optimum experimental design in this problem.
The last two cases (4 and 5) reflect an experiment with homogeneous Neumann boundary conditions ($\beta \equiv 0$). 
Comparing the results of case~1/2 and case~4, we see that the experiment of case~4 contains less information.
This can be attributed to the saturation in the diffusion process, which now lacks outflow.
The saturation eventually leads to equal concentrations in $\Omega$ and this gives rise to reduced sensitivities of the measurements with respect to interface perturbations on $\Gammainc$ at later time steps. 
As expected the non-optimized version of the Neumann boundary experiment with uniform measurement weights (case~5) yields the least amount of information among all experiments considered.

\begin{table}[htbp]
	\label{tab:compare_2D_settings}
	\begin{center}
		\begin{tabular}{ll}
			\toprule
			case~1:& optimized weights for $\beta=10$ on the lower left part of $\GammaR$\\
			case~2:& optimized weights for $\beta=10$ on the lower right part of $\GammaR$\\
			case~3:& uniform weights for $\beta=10$ on the lower left part of $\GammaR$\\
			case~4:& optimized weights with $\beta=0$ on $\GammaR$\\
			case~5:& uniform weights with $\beta=0$ on $\GammaR$\\
			\bottomrule
		\end{tabular}
		\caption{Overview of different experiments in the 2D case.}
	\end{center}
\end{table}

\begin{table}[htbp]
	\label{tab:compare_results_2D_settings}
	\begin{center}
		\begin{tabular}{c|c|c|c|c|c|c|c|c|c|c}
			case &
			$\scriptstyle{\Phi_A \cdot 10^{-4}}$ &
			\multicolumn{9}{c}{eigenvalues $(\Lambda_1^{-1},\ldots,\Lambda_9^{-1}) \cdot 10^{-4}$} \\
			\midrule
			1 &
			$32.93$ &
			$0.27$ &
			$0.49$ &
			$0.77$ &
			$1.04$ &
			$1.24$ &
			$2.58$ &
			$3.39$ & 
			$6.68$ &
			$16.48$ \\
			2 &
			$36.75$ &
			$0.39$ &
			$0.47$ & 
			$0.76$ & 
			$1.13$ & 
			$1.32$ & 
			$3.83$ &
			$5.53$ &
			$6.34$ &
			$16.97$ \\
			3 &
			$122.14$ &
			$0.57$ &
			$0.64$ & 
			$1.28$ & 
			$2.02$ & 
			$2.44$ & 
			$4.96$ &
			$10.31$ &
			$16.03$ &
			$83.89$ \\
			4 &
			$48.21$ &
			$0.42$ &
			$0.60$ & 
			$0.97$ & 
			$1.35$ & 
			$1.98$ & 
			$4.52$ &
			$6.85$ &
			$9.62$ &
			$21.89$ \\
			5 &
			$252.17$ &
			$1.36$ &
			$2.10$ & 
			$3.47$ & 
			$5.51$ & 
			$8.55$ & 
			$22.99$ &
			$30.78$ &
			$45.46$ &
			$131.96$ \\
			\bottomrule
		\end{tabular}
		\caption{Compare A-criterion \eqref{eq:A-criterion} and eigenvalues for different experiments, cf.\ \cref{tab:compare_2D_settings} for different cases.}
	\end{center}
\end{table}

\subsection{The 3D Case}
\label{subsec:Numerical_results_3D}

In three dimensions we consider a slightly modified test case.
The hold-all is given by $D = (0.1, 0.9)^3$ as a subregion of $\Omega = (0,1)^3$.
Further, in contrast to the 2D case, sensors are now flat and located on the outer boundary of $\Omega$.
Consequently, the distributed measurements in \eqref{eq:individual_measurement} need to be replaced by boundary measurements, i.e.,
\begin{equation}
	\label{eq:individual_measurements_3d}
	E_{k,\ell} \, u 
	:=
	\restr{u(t^\ell)}{\Gammaobs^k}
	\in 
	L^2(\Gammaobs^k)
\end{equation}
in the sense of traces.
Like in the 2D case we set $\alpha_0 = 0.01$ and $\alpha_1=1.0$ in \eqref{eq:covariance_root_differential_operator} to calculate $\AA_k$.
Notice that the covariance operator $\CC_k=\AA_k^{-2}$ is now defined on $\Gammaobs^k$.
Clearly, the inner product in \eqref{eq:elementary_FI_operator} needs to be replaced by the inner product in $L^2(\Gammaobs^k)$.

\begin{figure}[htbp]
	\label{fig:numbering_sensors_3d}
	\begin{center}
		\def\svgwidth{0.7\textwidth}
		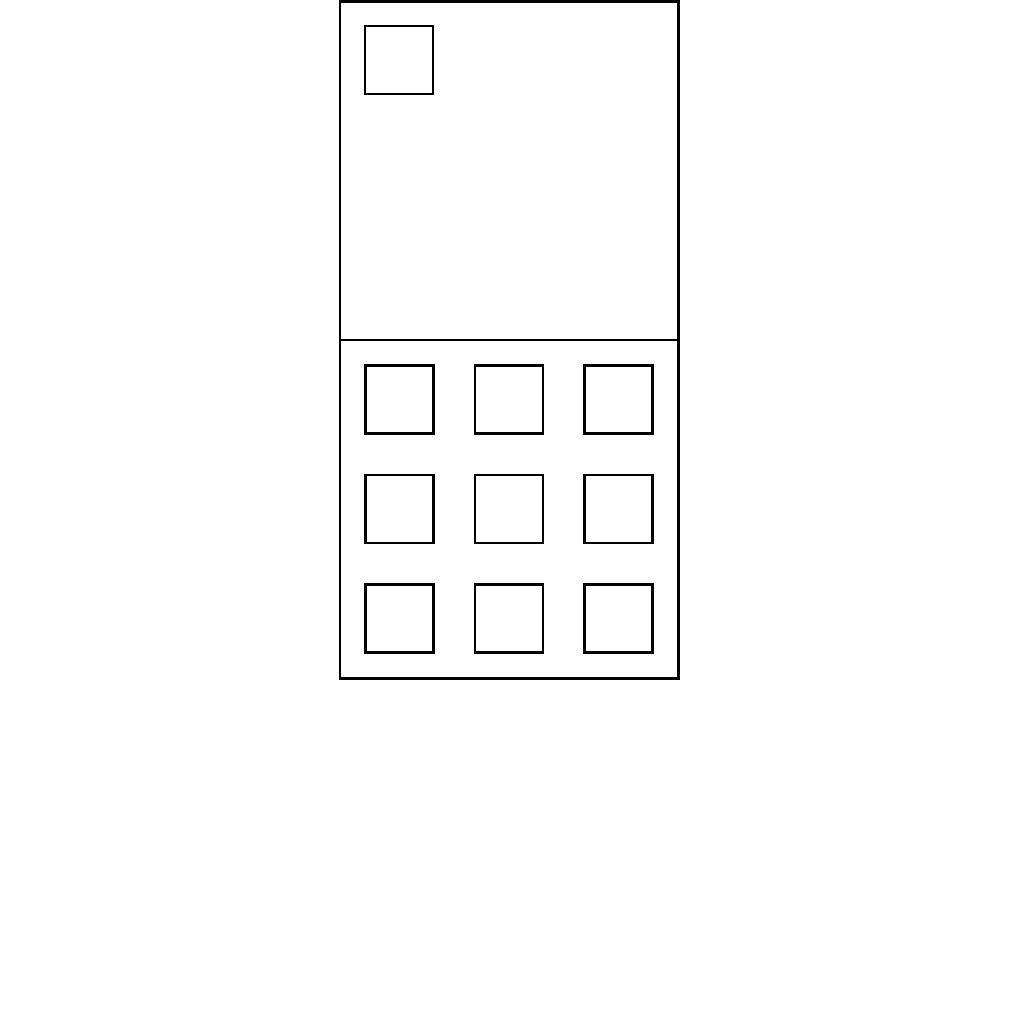
		\caption{Numbering of measurement sensors in the 3D experiment (unfolded unit cube $\Omega$), cf.\ \cref{fig:shape_sensors_3d}.}
	\end{center}
\end{figure}

On each facet of $\Omega$, except for the top, there are 9 quadratically shaped sensors with an edge length of $0.25$, which amounts to $\Nobs = 45$~sensors in total; see \cref{fig:numbering_sensors_3d}.
As the basis in the subspace of representative shape variations under consideration we use again deformation fields according to \eqref{eq:gaussian_bump}.
Here we use a slope factor of $s = 40$.
However, we have to find a different access to geodesic distances since we do not have a parametrization available for the triangulated surface mesh.
To this end we apply the Floyd--Warshall algorithm to the graph representing the surface triangulation of $\Gammainc$ to obtain pair-wise shortest paths between all finite element nodes; see \cite{Floyd1962,Hougardy2010}.
The result is a matrix with entries $\dist(x^i, x^j)$ for all surface finite element nodes $x^i$ and $x^j$ on $\Gammainc$.
For increasingly fine and sufficiently regular meshes, $\dist(x^i, x^j)$ approximates geodesic distances, which are then plugged into \eqref{eq:gaussian_bump}.
Compared to the 2D~case, it is not straightforward to choose equidistant center points for the Gaussian.
Instead, we use the following strategy to maximize pair-wise distances of the centers of the Gaussians in order to obtain a homogeneous distribution:
let $S$ denote the set of center points for Gaussians on the triangulated surface.
Further, let $X \subset \R^3$ denote the set of finite element nodes forming the surface.
Then \cref{alg:3d_basis_functions} generates homogeneously distributed shape variations.

\begin{algorithm}
	\begin{algorithmic}
		\State{$S \gets \lbrace x^0\rbrace$}
		\For{$j = 1,\dots, \Nbasis-1$}
		\State{$y \gets \underset{x \in X\setminus S}{\mathrm{argmax}} \; \min\limits_{\substack{x^i, x^j \in S \cup \lbrace x \rbrace\\ x^i \neq x^j}}\; \dist(x^i, x^j)$}
		\State{$S \gets S \cup \lbrace y \rbrace$}
		\EndFor
	\end{algorithmic}
	\caption{Finding homogeneously distributed Gaussian shape variations centered in the points $S$ on the triangulated surface $\Gammainc$.}
	\label{alg:3d_basis_functions}
\end{algorithm}

We choose $\Nbasis = 17$ as the dimension of the subspace of shape variations. The 17th shape variation is chosen to be the uniform normal vector field, representing a scaling of the shape, similar to the bottom right picture in \cref{fig:untransformed_bumps}.

\begin{figure}[htbp]
	\label{fig:shape_sensors_3d}
	\begin{center}
		\includegraphics[width=0.8\textwidth]{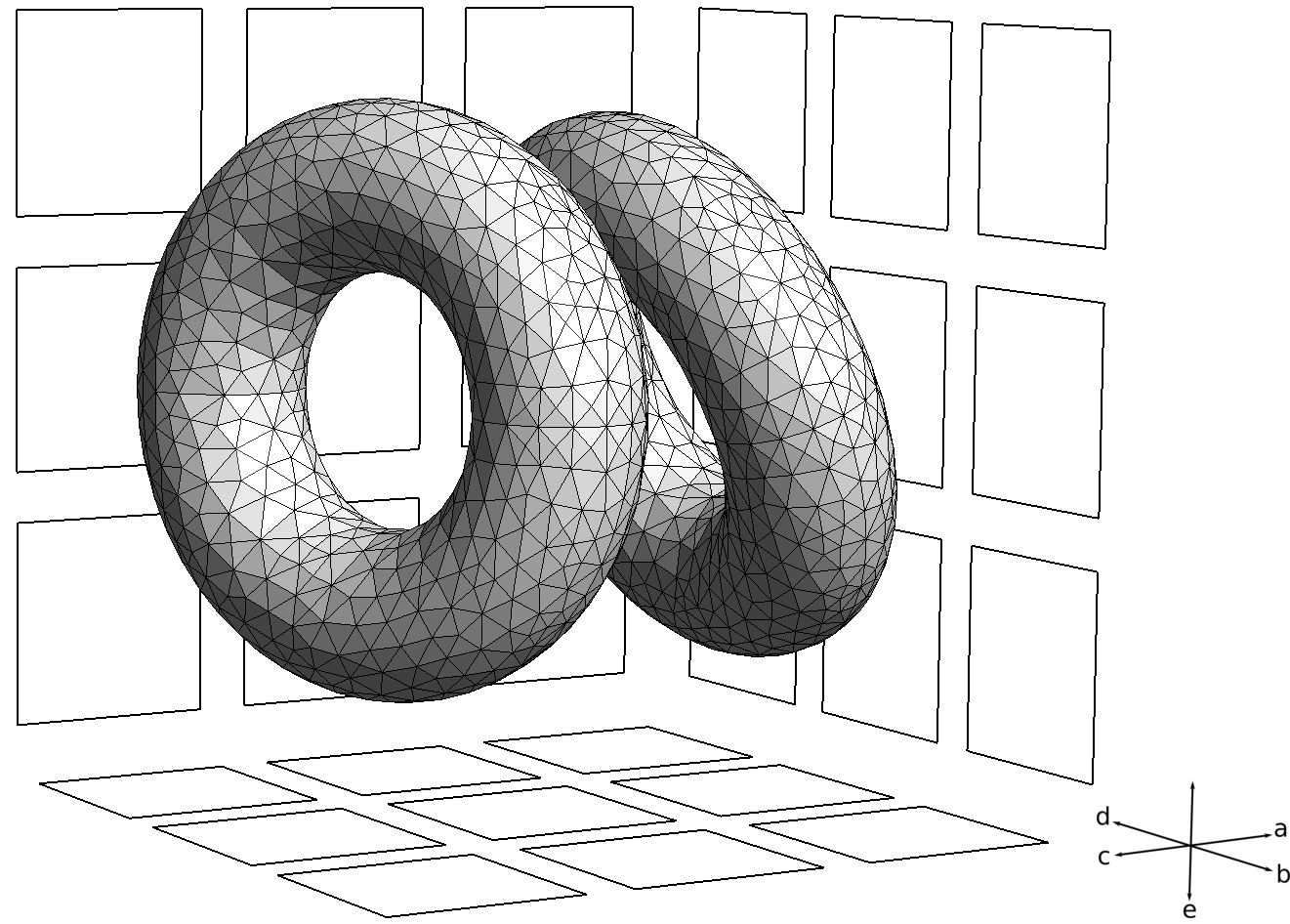}
		\caption{Interface $\Gammainc$ in 3D with some of the boundary measurement sensors, cf.\ \cref{fig:numbering_sensors_3d}.}
	\end{center}
\end{figure}

The underlying diffusion model \eqref{eq:forward_problem} is the same as in  the 2D setting with modified data.
The boundary conditions are now given by $\beta = 10$ on the entire bottom part of the boundary and $\beta = 0$ on the remaining parts of $\GammaR$.
The final time is chosen to be $T = 7$, and the interval $[0,T]$ is discretized into 21~equal intervals.
Again all time steps are assumed to be potential points for measurements yielding $\Ntime = 22$.
Consequently, we have $\Nobs \cdot \Ntime = 990$~potential sensor activations.
We choose the upper bound for the sum of the weights in the optimum experimental design to be $C_w = 40$.

\begin{figure}[htbp]
	\begin{center}
		\includegraphics[width=0.48\textwidth]{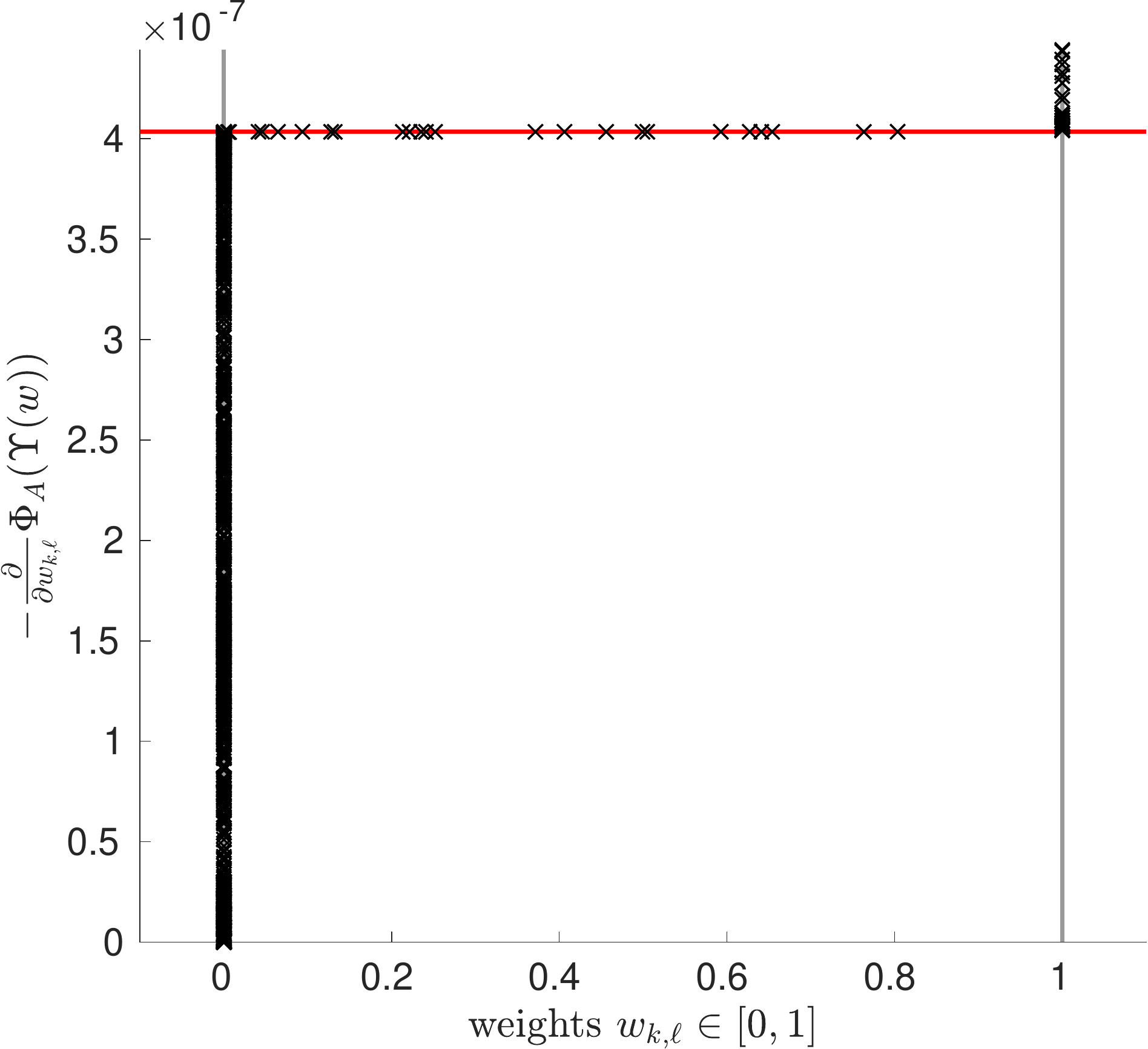}
		\hfill
		\includegraphics[width=0.48\linewidth]{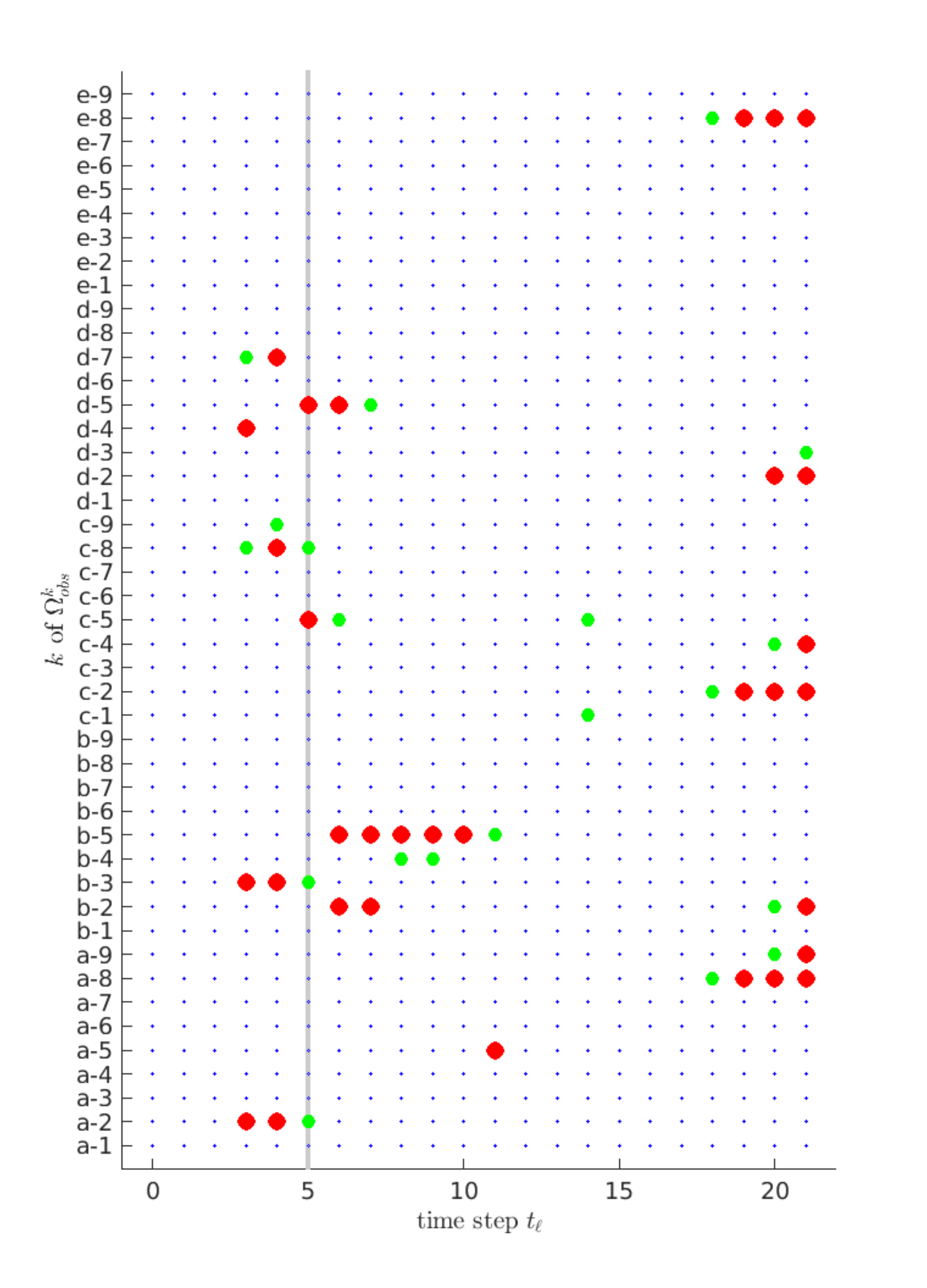}
	\end{center}
	\caption{Results for the 3D experiment: the relaxed optimality condition \eqref{eq:optimality_condition_for_relaxed_OED_problem} is verified at the final iterate (left); optimal sensor activation pattern (right) with $w_{k,\ell}=1$ shown in red and $w_{k,\ell} \in \left(0,1\right)$ shown in green. The sensor activation pattern at the 5th~time step (indicated by the gray vertical line) is shown in \cref{fig:sensor_activation_t5}. The measurement domains $\Omegaobs^k$ are numbered as in \cref{fig:numbering_sensors_3d}.}
	\label{fig:optimality_condition_3D}
\end{figure}

The optimization algorithm is the same as in the 2D case.
Again we choose a homogeneous distribution of the total weight $C_w$ as an initial guess.
After 154 iterations of the simplicial decomposition algorithm, an inexact version of the optimality condition \eqref{eq:optimality_condition_for_relaxed_OED_problem} is fulfilled.
This is visualized in \cref{fig:optimality_condition_3D}, where the red horizontal line marks the value of~$\xi$ in \eqref{eq:optimality_condition_for_relaxed_OED_problem}.
Like in the 2D case, the optimized weight vector $w^*$ exhibits a sparse structure, which is shown in \cref{fig:optimality_condition_3D}.
In total we encounter 32~weights equal to one, 20 in $(0,1)$ and 938 are zero.
In the final iteration of the simplicial decomposition algorithm there are 23~active vertices spanning a subset of $\Delta_{C_w}$ (cf.\ \cref{fig:restricted_simplex} and \eqref{eq:simplex}), which includes $w^\ast$.

\begin{figure}[htbp]
	\label{fig:optimized_bumps_3d}
	\begin{center}
		\begin{tabular}{ccc}
			\includegraphics[width=.4\textwidth]{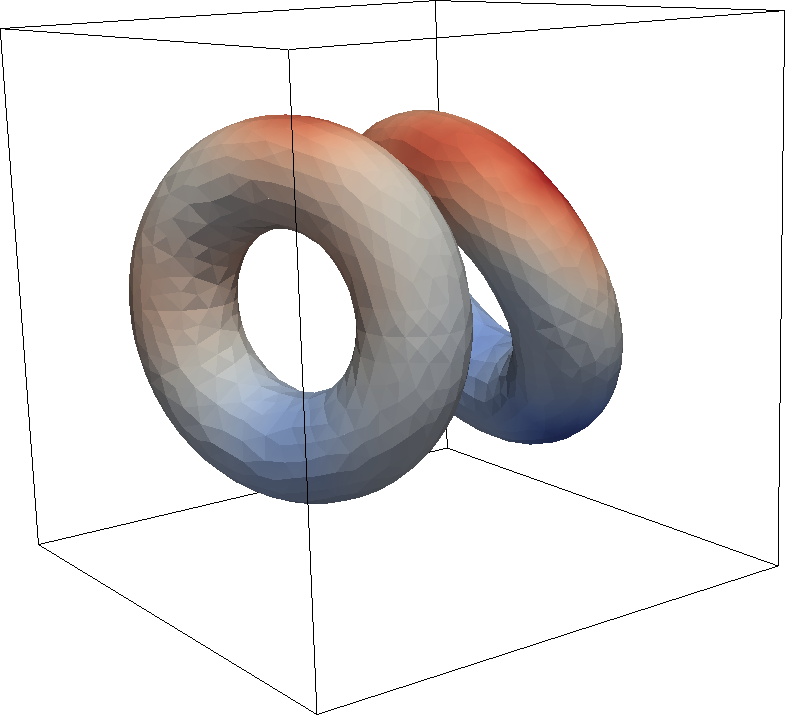}&
			\includegraphics[width=.4\textwidth]{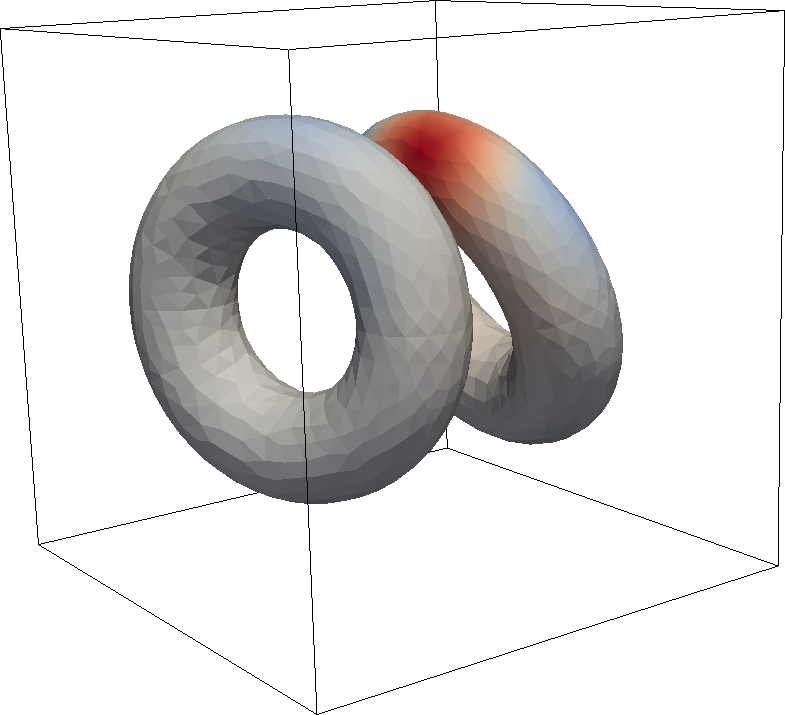}&
			\includegraphics[width=.06\textwidth]{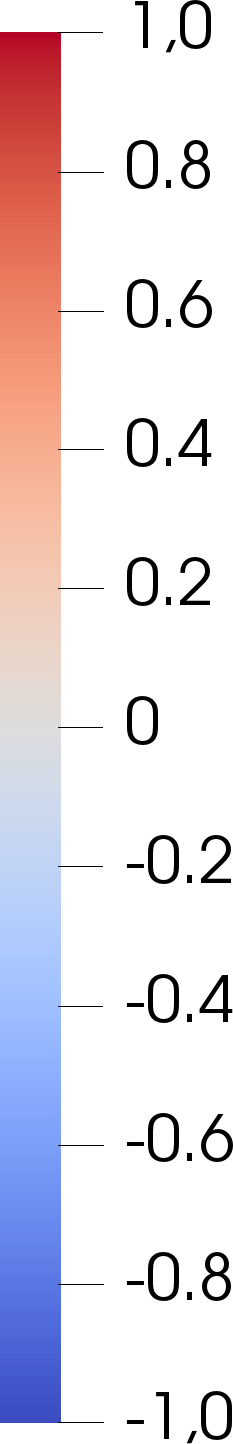}\\
		\end{tabular}
		\caption{Eigenfunctions (orthonormal with respect to $B$ from \eqref{eq:generalized_eigenvalue_problem}, but shown here with unified scale) for the smallest (left) and largest (right) inverse eigenvalue according to the optimized weights. The vector field defined by the respective eigenfunction is normal to the surface $\Gammainc$. A positive/negative value (red/blue) indicates an outward/inward pointing normal vector. The corresponding eigenvalues $(\Lambda_1^{-1},\Lambda_{17}^{-1})$ are $(0.20,49.34)\times\num{e-7}$.}
	\end{center}
\end{figure}

In order to geometrically investigate the identifiability of the inclusion $\Gammainc$ we solve the generalized eigenvalue problem \eqref{eq:generalized_eigenvalue_problem} for the FIM $\Upsilon(w^\ast)$ evaluated at the optimized weight vector $w^\ast$.
In \cref{fig:optimized_bumps_3d} we can see two eigenfunctions where the color encodes the normal component. 
On the left hand side we show the eigenfunction corresponding to the smallest reciprocal eigenvalue.
It can be interpreted as the direction in the subspace of the deformation tangent space corresponding to the best identifiability.
On the right hand side, by contrast, we display the eigenfunction of worst identifiability.

\begin{figure}[htbp]
	\label{fig:optimal_eigenvalues_3d}
	\begin{center}
		{%
			\setlength{\fboxsep}{7pt}%
			\setlength{\fboxrule}{1pt}%
			\fbox{\includegraphics[width=.6\textwidth]{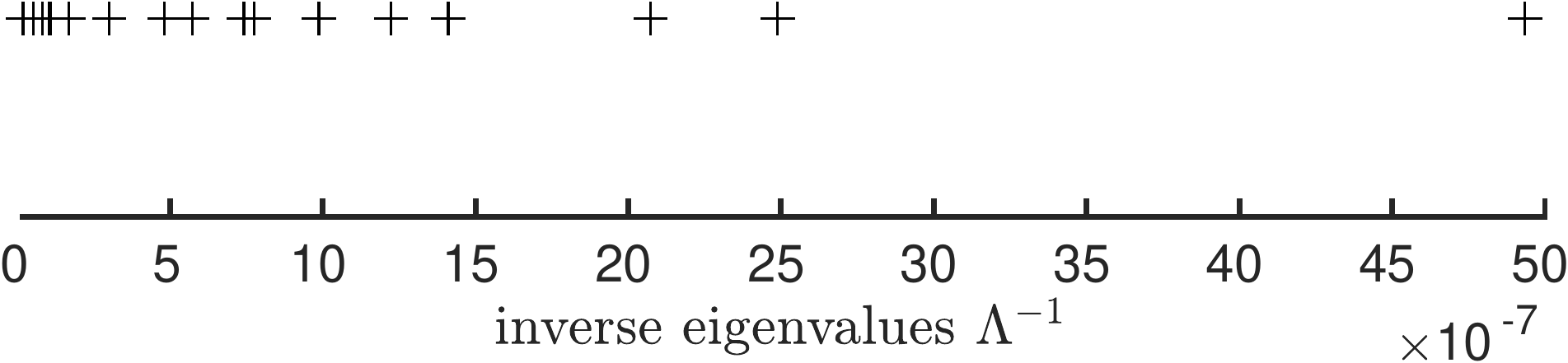}}
		}
		\caption{Inverse eigenvalues $(\Lambda_1^{-1},\ldots,\Lambda_{17}^{-1})$ according to the FIM of the optimized weights in the 3D experiment. 
			The smallest and largest eigenvalues are $\num{0.20e-7}$ and $\num{49.34e-7}$.
		}
	\end{center}
\end{figure}

We can interpret these results geometrically as follows.
A redistribution of volume between the top and bottom halves of $\Omegainc$ is well identifiable by the optimized experiment.
By contrast, localized shape modifications near the top, where the distance between the tori is minimal and the distance to the sensors is maximal, are hard to identify.
\Cref{fig:optimal_eigenvalues_3d} visualizes the corresponding eigenvalue distribution.
Like in the 2D case, we can observe a significant gap between the largest reciprocal eigenvalue and the next better ones.
This indicates that the eigenfunction in the right of \cref{fig:optimized_bumps_3d} represents by far the worst identifiable shape variation.

\begin{figure}[htbp]
	\begin{center}
		\includegraphics[width=0.7\linewidth]{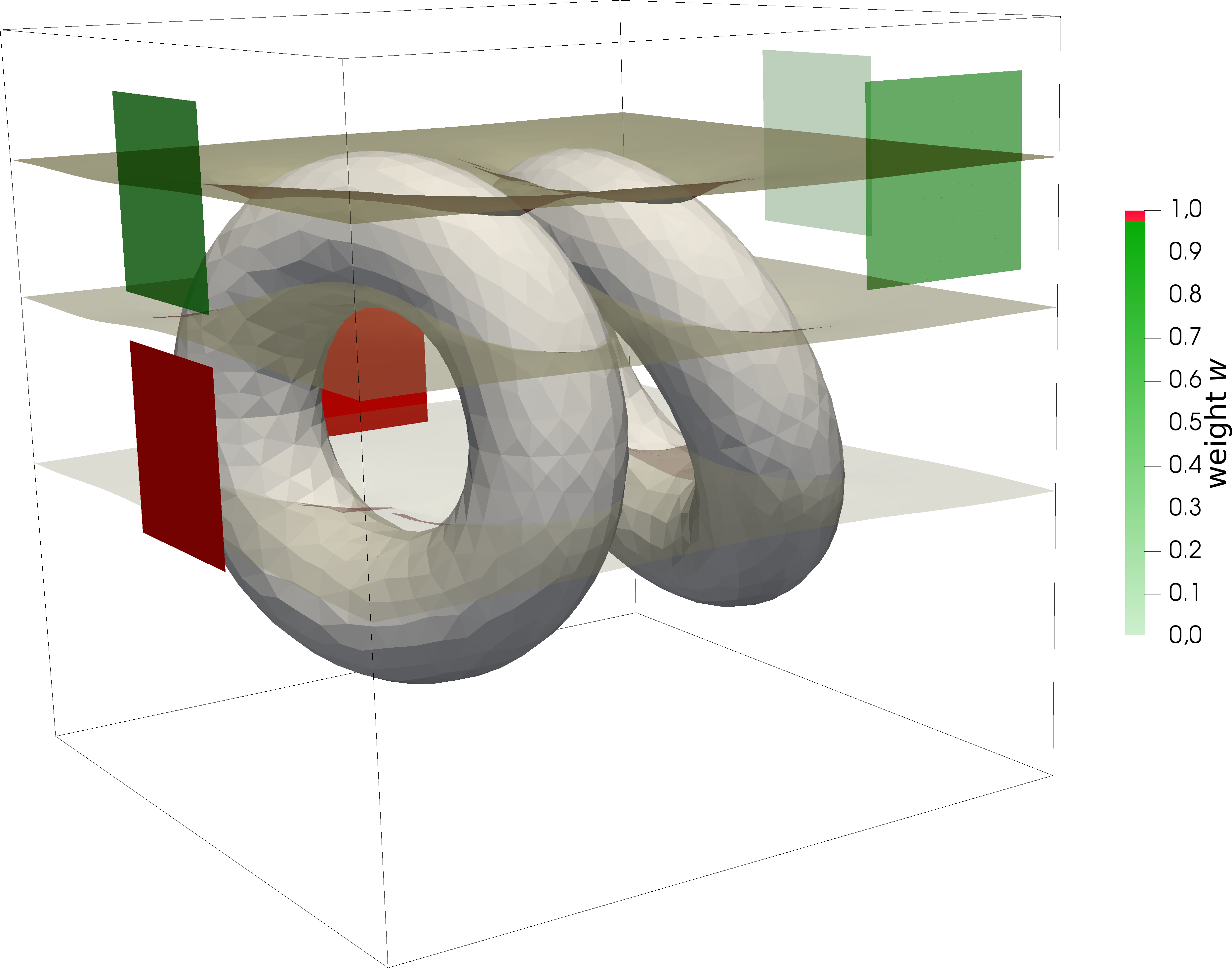}
	\end{center}
	\caption{Optimal sensor activation pattern at time step 5, with iso-surfaces at concentrations $\{750,500,250\}$. There are two activated sensors with $w_{k,\ell}=1$ at positions $\{c_5,d_5\}$ and three sensors with $w_{k,\ell} \in \left(0,1\right)$ at $\{a_2,b_3,c_8\}$. The numbering is as in \cref{fig:numbering_sensors_3d}.}
	\label{fig:sensor_activation_t5}
\end{figure}

\subsection{Spatial-Only Sensor Activation (3D)}
\label{subsec:Numerical_results_spatial}

In this third part of our numerical tests we demonstrate that the general setting of \cref{sec:OED_problem} can be easily modified to purely spatial sensor placement problems.
This means that we seek to select sensors which then remain activated over the entire time horizon $[0,T]$.
Owing to the independence of measurement outcomes at individual regions and time instances, the combined Fisher information matrix (FIM) \eqref{eq:combined_experiment} becomes
\begin{equation}
	\label{eq:combined_experiment_spatial}
	\Upsilon_\text{spatial}(w) \coloneqq \sum_{k=1}^{\Nobs} w_k \left[\sum_{\ell=1}^{\Ntime} \Upsilon_{k,\ell}\right].
\end{equation}
Note that no changes to the covariance operator $\CC_k$ have to be made since we only consider spatial and no temporal correlations of measurements.
We formulate the spatial-only OED problem as
\begin{equation}
	\label{eq:relaxed_OED_problem_spatial}
	\begin{aligned}
		\text{Minimize} \quad & \Phi_A (\Upsilon_\text{spatial}(w))\\
		\text{s.t.} \quad & 0 \le w_{k} \le 1 \quad \text{for all } k = 1, \ldots, \Nobs\\
		\text{and} \quad & \sum_{k=1}^{\Nobs} w_{k} \le C_w.
	\end{aligned}
\end{equation}
Problem \eqref{eq:relaxed_OED_problem_spatial} is readily solved by the same simplicial decomposition algorithm used for the experiments in \cref{subsec:Numerical_results_2D} and \cref{subsec:Numerical_results_3D}.
The only change necessary is to replace the space-time elementary FIMs $\Upsilon_{k,\ell}$ with their summation in time, i.e., $\Upsilon_k \coloneqq \sum_{\ell=1}^{\Ntime} \Upsilon_{k,\ell}$, $k=1,\dots,\Nobs$.
The dimension of the optimization problem reduces from $\Ntime\cdot \Nobs$ to $\Nobs = 45$.

\begin{figure}[htbp]
	\begin{center}
		\includegraphics[width=0.7\linewidth]{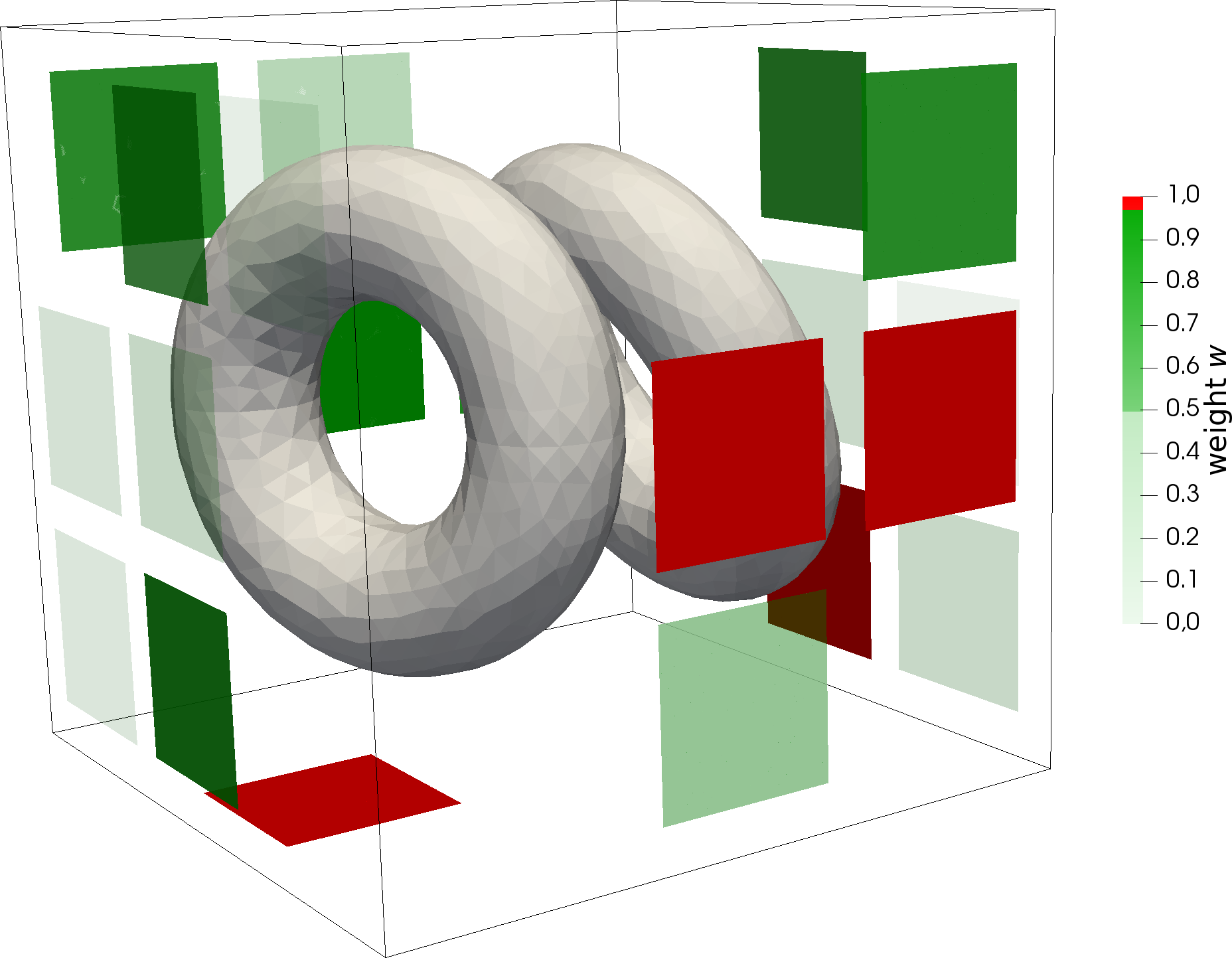}
	\end{center}
	\caption{Optimal spatial-only sensor activation pattern. There are four activated sensors with $w_{k}=1$ at $\{a_8,b_2,b_5,e_8\}$ and 16 sensors with $w_{k} \in \left(0,1\right)$. The numbering is as in \cref{fig:numbering_sensors_3d}.}
	\label{fig:spatial_sensor_activation}
\end{figure}

\begin{figure}[htbp]
	\label{fig:optimized_spatial_bumps_3d}
	\begin{center}
		\begin{tabular}{ccc}
			\includegraphics[width=.4\textwidth]{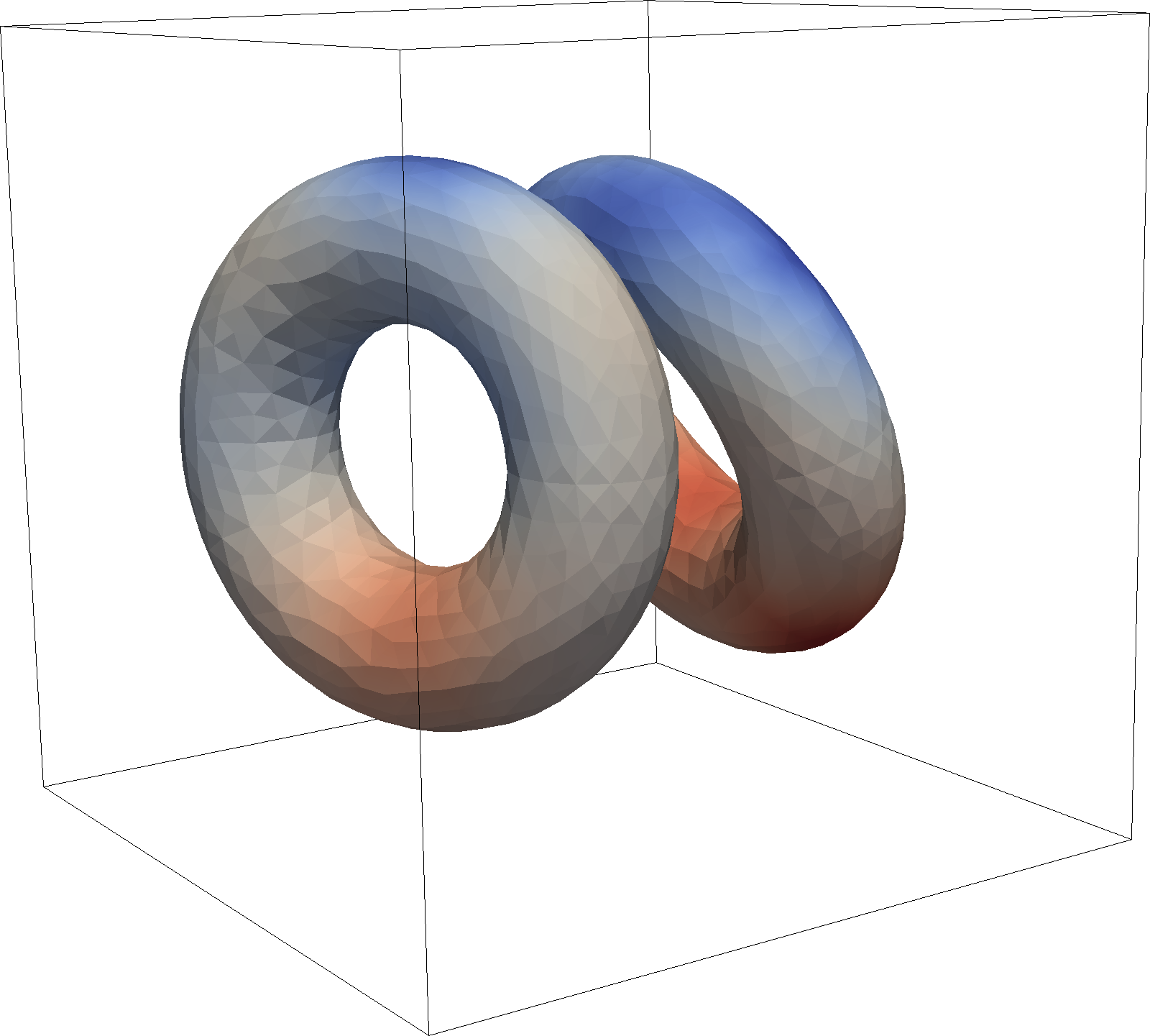}&
			\includegraphics[width=.4\textwidth]{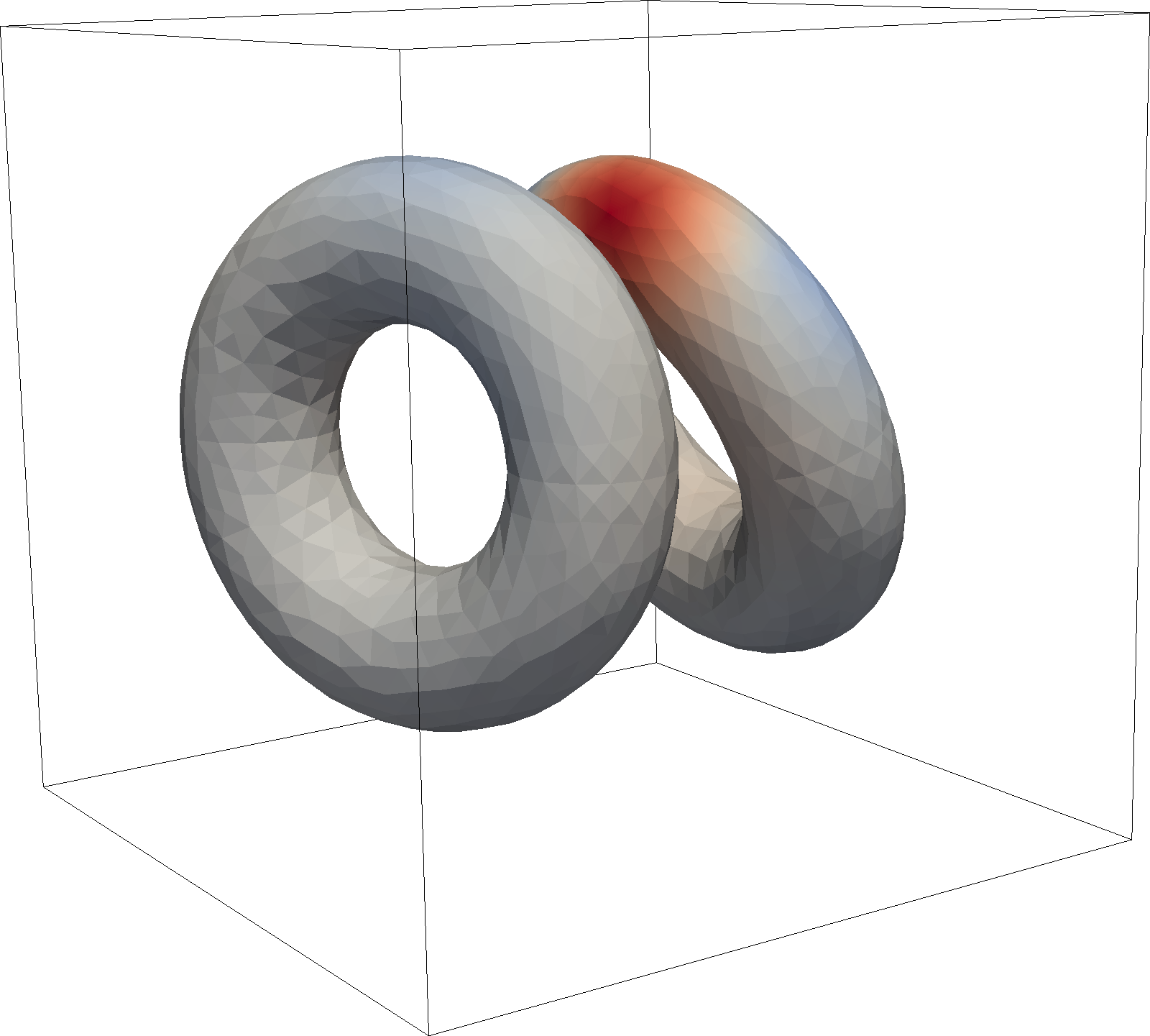}&
			\includegraphics[width=.06\textwidth]{Pictures/3d_bumps/colorbar_vert_ef_3d.png}\\
		\end{tabular}
		\caption{Eigenfunctions (orthonormal with respect to $B$ from \eqref{eq:generalized_eigenvalue_problem}, but shown here with unified scale) for the smallest (left) and largest (right) inverse eigenvalue according to the optimized weights. The vector field defined by the eigenfunctions points in the normal directions of the surface $\Gammainc$. A positive value (red) means that it points outwards and a negative value (blue) inwards. The corresponding eigenvalues $(\Lambda_1^{-1},\Lambda_{17}^{-1})$ are $(0.03,11.82)\times\num{e-7}$.}   
	\end{center}
\end{figure}

\begin{figure}[htbp]
	\label{fig:optimal_eigenvalues_3d_spatial}
	\begin{center}
		{
			\setlength{\fboxsep}{7pt}
			\setlength{\fboxrule}{1pt}
			\fbox{\includegraphics[width=.6\textwidth]{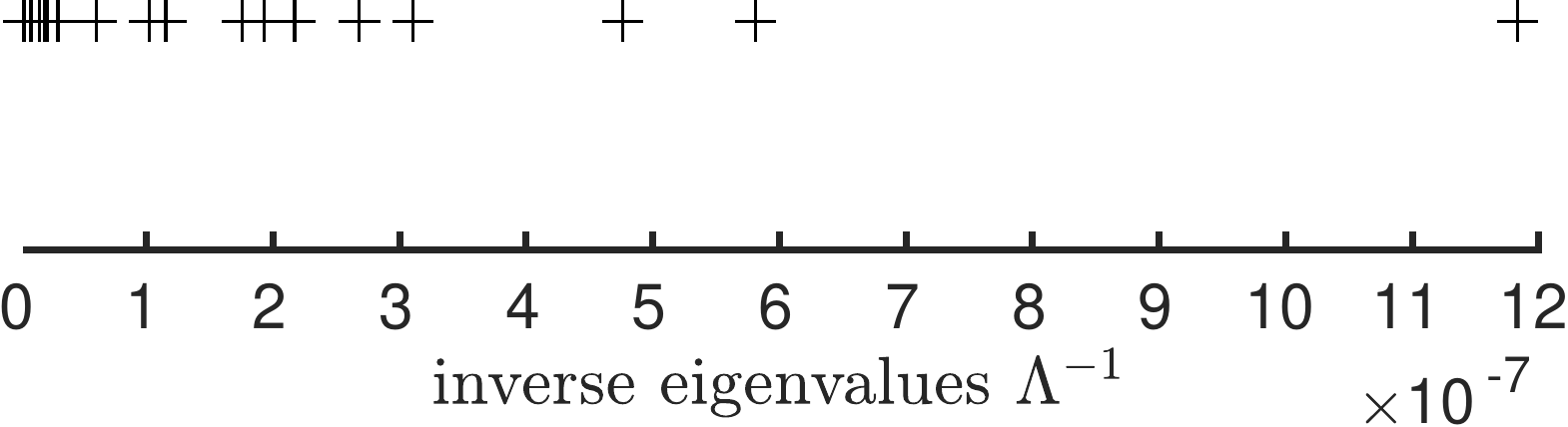}}
		}
		\caption{Inverse eigenvalues $(\Lambda_1^{-1},\ldots,\Lambda_{17}^{-1})$ according to the FIM of the optimized weights in the 3D spatial-only sensor activation experiment.
			The smallest and largest eigenvalues are $\num{0.03e-7}$ and $\num{11.82e-7}$, respectively.
		}
	\end{center}
\end{figure}

In our experiment we allow a maximum total weight of $C_w = 12$ whereas the other settings are identical to those in \cref{subsec:Numerical_results_3D}.
The optimal activation pattern is visualized in \cref{fig:spatial_sensor_activation}. 
There are four fully activated sensors with $w_{k}=1$ at $\{a_8,b_2,b_5,e_8\}$ and 16 sensors with $w_{k} \in \left(0,1\right)$.
The corresponding best and worst identifiable eigenfunctions are shown in \cref{fig:optimized_spatial_bumps_3d}.
Note that these are similar to the results obtained in the space-time sensor activation problem.
\Cref{fig:optimal_eigenvalues_3d_spatial} visualizes the inverse eigenvalues corresponding to the optimized FIM.
Their distribution is comparable to the results in \cref{subsec:Numerical_results_3D}.

\section{Discussion and Conclusion}
\label{sec:Discussion_conclusion}

In this work we have presented theory and an algorithm for optimum experimental design for interface estimation problems, which can be viewed as parameter estimation problem in which the parameter space does not exhibit a vector space but rather a manifold structure.
As a particular example, we considered optimal sensor activation problems for two and three dimensional interface identification experiments in a diffusion process.
A natural extension of this model problem is to combine it with additional experimental conditions such as the choice of the Robin parameter $\beta$ along the outflow boundary.
In this setting, the elementary Fisher information matrices depend on $\beta$ and thus they can no longer be pre-calculated.
The efficient solution of this extended problem is left to future research.


\IfFileExists{World.bib}
{
	\bibliographystyle{siamplain}
	\bibliography{World}
}
{
	\bibliographystyle{siamplain}
	\bibliography{paper}
}

\end{document}